\documentclass{article}

\usepackage{amsmath,amssymb, verbatim}
\usepackage{graphicx} 

\newtheorem{theorem}{Theorem}[section]
\newtheorem{corollary}[theorem]{Corollary}
\newtheorem{lemma}[theorem]{Lemma}
\newtheorem{proposition}[theorem]{Proposition}
\newtheorem{definition}[theorem]{Definition}

\newtheorem{remark}[theorem]{Remark}
\numberwithin{equation}{section}

\def\med{\medbreak\noindent}
\def\sms{\smallskip}
\def\sm{\smallskip\noindent}
\def\ms{\medskip}

\def\small{\smallskip\noindent}
\def\nl{\hfil\break}
\def\noi{\noindent}

\def\today{\noindent\number\day
\space\ifcase\month\or
  January\or February\or March\or April\or May\or June\or
  July\or August\or September\or October\or November\or December\fi
  \space\number\year}

 \def\bE {{\mathbb E}}

 \def\bN {{\mathbb N}} 
\def\bP {{\mathbb P}} \def\bQ {{\mathbb Q}} \def\bR {{\mathbb R}}

 \def\bZ {{\mathbb Z}}

\def\sA {\mathcal{A}} \def\sB {\mathcal{B}} \def\sC {\mathcal{C}}
\def\sD {\mathcal{D}} \def\sE {\mathcal{E}} \def\sF {\mathcal{F}}
 \def\sH {\mathcal{H}} 
  \def\sL {\mathcal{L}}
 \def\sN {\mathcal{N}} 
 \def\sQ {\mathcal{Q}} \def\sR {\mathcal{R}}
\def\sS {\mathcal{S}}

\def\ol{\overline}

\def\al {\alpha}
\def\lam {\lambda} 

\def\gam{\gamma}
\def\Gam{\Gamma}

\def\eps{\varepsilon}
\def\om{\omega }
\def\th{\theta} 
\def\Th{\Theta}
\def\vp{\varphi}
\def\phi{\varphi}

\def\pd {\partial}
\def\q{\quad} \def\qq{\qquad}
\def\dint{\int\kern-.6em\int}
\def\grad{\nabla}

\def\inter{\mathop{{\rm int }}}

\def\supp{\mathop{{\rm supp}}}

\def\Cap{\mathop{{\rm Cap }}}

\def\inter{{\mathop {{\rm int}}}}

\def \half {{\textstyle\frac12}}

\def\=d{{\,\buildrel (d) \over =\,}}
\def\a.s.{{\buildrel a.s. \over \longrightarrow}}
  
\def\tilm{{\sim_m}}

\def\wt{\widetilde}

\def\wbP{\wt \bP}

\def\q{\quad}
\def\qq{{\qquad}}
\def\nl{{\newline}}
\def\pd{\partial}

\def\proof{{\medskip\noindent {\bf Proof. }}}

\def\square{{\vcenter{\vbox{\hrule height.3pt
                  \hbox{\vrule width.3pt height5pt \kern5pt
                     \vrule width.3pt}
                  \hrule height.3pt}}}}

\def\qed{{\hfill $\square$ \medskip}}

\def\nn{\nonumber}

\def\fE{{\mathfrak E}}
\def\Reff{R_{\rm eff}}
\def\sEBB{\sE_{BB}}
\def\sEKZ{\sE_{KZ}}
\def\LF{L_F}
\def\MF{m_F}
\def\h{{h}}

\begin{document}

\title{\bf Uniqueness of Brownian motion on Sierpinski carpets}

\author{
{\bf Martin T. Barlow}
\thanks{Research partially supported by NSERC (Canada), and EPSRC (UK).}
\\
Department of Mathematics, University of British Columbia\\ 
Vancouver B.C. Canada V6T 1Z2\\ 
Email: barlow@math.ubc.ca \medskip \\
{\bf Richard F. Bass}
\thanks{Research partially supported by NSF grant DMS-0601783.}
\\
Department of Mathematics, University of Connecticut\\
Storrs~CT~06269-3009~USA\\ 
Email: bass@math.uconn.edu \medskip \\ 
{\bf Takashi Kumagai}
\thanks{Corresponding author}
\thanks{Research partially supported by the Grant-in-Aid
for Scientific Research (B) 18340027.}\\
Department of Mathematics, Faculty of Science\\
Kyoto University, Kyoto 606-8502, Japan\\ 
Email: kumagai@math.kyoto-u.ac.jp \medskip \\ 
and \medskip \\
{\bf Alexander Teplyaev}\thanks{Research partially supported by  NSF grant
DMS-0505622.}\\
Department of Mathematics, University of Connecticut\\
Storrs~CT~06269-3009~USA\\
Email: teplyaev@math.uconn.edu}

\maketitle

\begin{abstract}
We prove that, up to scalar multiples, there exists only one 
local regular Dirichlet
form on a generalized Sierpinski carpet that is invariant with respect
to the local symmetries of the carpet.  Consequently for each such
fractal the law of Brownian motion is uniquely determined and the 
Laplacian is well defined.
%
%
\end{abstract}

\section{Introduction}

The standard Sierpinski carpet $F_{\rm SC}$ is the  fractal that is formed by
taking the unit square, dividing it into 9 equal subsquares,
removing the central square, dividing each of the 8 remaining
subsquares into 9 equal smaller pieces, and continuing. 
In \cite{BB1} two of the authors of this paper 
gave a construction of  a Brownian motion on $F_{\rm SC}$. 
This is a diffusion (that is, a continuous strong Markov process) 
which takes its values in  $F_{\rm SC}$, and which is 
non-degenerate and invariant under all the local isometries of 
$F_{\rm SC}$.

Subsequently, Kusuoka and Zhou in \cite{KZ} gave a different construction
of a diffusion 
on $F_{\rm SC}$, which yielded a process that, as well as having the
invariance properties of the Brownian motion constructed in \cite{BB1}, 
was also scale invariant.
The proofs in \cite{BB1, KZ} also work for fractals that are formed in a 
similar manner to the standard Sierpinski carpet: we call these 
{\sl generalized Sierpinski carpets} 
(GSCs). In \cite{BB4} the results of \cite{BB1} were 
extended to GSCs embedded in $\bR^d$ for $d\geq 3$.
While \cite{BB1,BB4} and \cite{KZ} both obtained their
diffusions as limits of
approximating processes, the type of approximation was different: 
\cite{BB1,BB4} used a sequence of time changed reflecting Brownian motions, 
while \cite{KZ} used a sequence of Markov chains. 

\begin{figure} [h]
\centering
  \includegraphics[width=60mm] {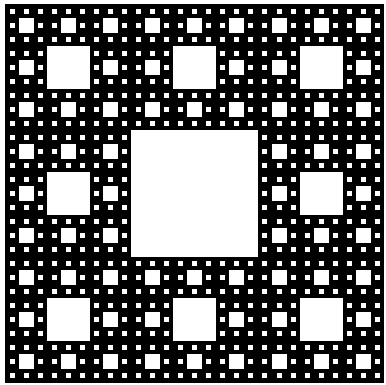}
\caption{ The standard Sierpinski carpet}
\end{figure}

These papers left open the question
of uniqueness of this Brownian motion -- in fact it was not even 
clear whether or not the processes obtained in \cite{BB1, BB4} or  \cite{KZ}
were the same. 
This uniqueness question can also be expressed in analytic terms: one can
define a {\sl Laplacian} on a GSC as the infinitesimal
generator of a Brownian motion, and one wants to know if 
there is only one such Laplacian. The main result of this paper is that, up
to scalar multiples of the time parameter,
there exists only one such Brownian motion; hence,
up to scalar multiples, the Laplacian is uniquely defined. 

GSCs are examples of spaces with {\sl anomalous diffusion}. For Brownian 
motion on $\bR^d$ one has $\bE |X_t-X_0| = c t^{1/2}$. 
Anomalous diffusion in a space $F$ occurs when instead one has 
$\bE |X_t-X_0| = o(t^{1/2})$, or (in regular enough situations),
$\bE |X_t-X_0| \approx t^{1/d_w}$, 
where $d_w$ (called the {\sl walk dimension})
satisfies $d_w>2$. This phenomena was first observed by mathematical 
physicists working in the transport properties of disordered media, such as 
(critical) percolation clusters -- see \cite{AO, RT}. Since these sets are 
subsets of the lattice $\bZ^d$, they are not true fractals, but their large scale
structure still exhibits fractal properties, and the simple random walk
is expected to have anomalous diffusion.

For critical percolation clusters (or, more precisely for the
incipient infinite cluster) on trees and $\bZ^2$, Kesten \cite{Kes1} proved 
that anomalous diffusion occurs. After this work,
little progress was made on critical percolation clusters 
until the recent papers \cite{BK, BJKS, KN}.

As random sets are hard to study, it was natural to begin the study 
of anomalous diffusion in the more tractable context of 
regular deterministic fractals. The simplest of these is the Sierpinski gasket.
The papers \cite{AO, RT} studied discrete random 
walks on graph approximations to the Sierpinski gasket, and soon after \cite{G, Kus0, BP}
constructed Brownian motions on the limiting set. 
The special structure of the Sierpinski gasket makes the uniqueness problem
quite simple, and uniqueness of this Brownian motion was proved in \cite{BP}.
These early papers used a probabilistic approach, first constructing
the Brownian motion $X$ on the space, and then, having defined the Laplacian 
$\sL_X$ as the infinitesimal generator of the semigroup of $X$, used the process
$X$ to study  $\sL_X$. Soon after Kigami \cite{Kigsg} and 
Fukushima-Shima \cite{FS} introduced more analytical approaches, and
in particular \cite{FS} gave a very simple construction of $X$ and $\sL_X$
using the theory of Dirichlet forms. 

\begin{figure} [h]
\centering
  \includegraphics[width=120 mm]{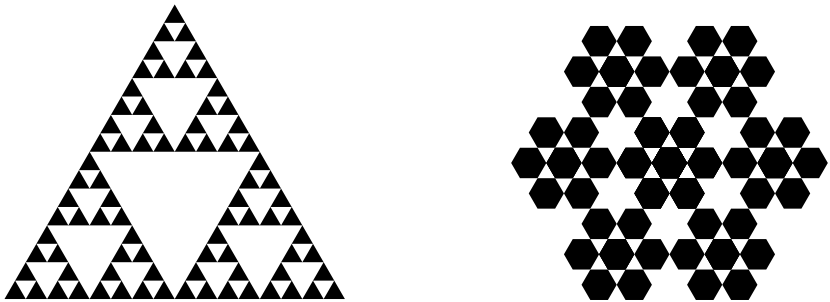}
\caption{ The Sierpinski gasket (left), and a typical nested fractal, the
Lindstr\o m snowflake (right)} 
\end{figure}

It was natural to ask whether these results were special to the Sierpinski gasket.
Lindstr\o m \cite{L} and Kigami \cite{Kigpcf} introduced wider families of
fractals (called {\sl nested fractals}, and {\sl p.c.f.~self-similar sets} respectively), and
gave constructions of diffusions on these spaces. Nested fractals are, like the
Sierpinski carpet, highly symmetric, and the uniqueness problem can be
formulated in a similar fashion to that for GSCs. Uniqueness for nested fractals was 
not treated in \cite{L}, and for some years remained 
a significant challenge, before being solved by Sabot \cite{Sab1}. (See also
\cite{Me1, Peirone2000}  for shorter proofs). 
For p.c.f.~self-similar sets, while some sufficient conditions for
uniqueness are given in \cite{Sab1, HMT}, the general problem is still open. 

The study of these various families of fractals (nested fractals, p.c.f self-similar
sets, and GSCs) revealed a number of common themes, and showed that
analysis on these spaces differs from that in standard
Euclidean space in several ways, all ultimately connected with the fact that
$d_w>2$:
\begin{itemize}
\item  The energy measure $\nu$
and the Hausdorff measure $\mu$ are mutually singular,
\item The domain of the Laplacian is not an algebra, 
\item If $d(x,y)$ is the shortest path metric, then $d(x, \cdot)$ is not
in the domain of the  Dirichlet form.
\end{itemize}
See \cite{bar, kig, Stri} for further information and references.

The uniqueness proofs in \cite{HMT, Me1, Peirone2000, Sab1} all used in an
essential way the fact that 
nested fractals and p.c.f.~self-similar sets are finitely ramified --
that is, they can be disconnected by removing a finite number of points.
For these sets there is a natural
definition of a set $V_n$ of `boundary points at level $n$' -- for the Sierpinski gasket
$V_n$ is the set of vertices of triangles of side $2^{-n}$. 
If one just looks at the process $X$ at the times when it passes through the
points in $V_n$, one sees a finite state Markov chain $X^{(n)}$, which
is called the {\em trace of $X$ on $V_n$}. 
If $m>n$ then $V_n \subset V_m$ and 
the trace of $X^{(m)}$ on $V_n$ is also $X^{(n)}$.
Using this, and the fact that the limiting processes are known to be scale 
invariant, the uniqueness problem for $X$ can be reduced to the 
uniqueness of the fixed point of a non-linear map on a space
of finite matrices. 

While the boundaries of the squares (or cubes) 
have an analogous role to the sets $V_n$ in the geometrical construction 
of a GSC, attempts to follow the same strategy of proof 
encounter numerous difficulties and have
not been successful.
We use a different idea in this paper, and rather than studying the
restriction of the process $X$ to boundaries, our argument
treats the Dirichlet form of the process on the whole space.
(This also suggests a new approach to uniqueness on finitely ramified
fractals, which will be explored elsewhere.)

\ms Let $F$ be a GSC and $\mu$ the usual Hausdorff measure on $F$.
Let $\fE$ be the set of non-zero local regular conservative 
Dirichlet forms $(\sE, \sF)$ on $L^2(F,\mu)$ which are invariant with
respect to all the local symmetries of $F$.  
(See Definition \ref{deffE} for a precise definition.) We remark 
 that elements of
$\fE$ are not required to be scale invariant -- see Definition~\ref{def-scaleinv}.
Our first result is that $\fE$ is non-empty.

\begin{proposition} \label{ebbkz}
The Dirichlet forms associated with the processes 
constructed in \cite{BB1, BB4} and \cite{KZ} are in $\fE$.
\end{proposition}

Our main result is the following theorem, which is proved in
Section~\ref{Uni}.

\begin{theorem} \label{tmain}
Let $F\subset \bR^d$ be a GSC.
Then, up to scalar multiples, $\fE$ consists of at most one element.
Further, this one element of $\fE$ satisfies scale invariance.
\end{theorem}

An immediate corollary of Proposition \ref{ebbkz} 
and Theorem \ref{tmain} is the following.

\begin{corollary}\label{C1.2}
 The Dirichlet forms constructed in
\cite{BB1, BB4} and \cite{KZ} are (up to a constant) the same.
\newline (b) The Dirichlet forms constructed in \cite{BB1, BB4}
 satisfy scale invariance.
\end{corollary}

A Feller process is one where the semigroup $T_t$ maps continuous
functions that vanish at infinity to continuous functions that vanish
at infinity, and $\lim_{t\to 0}T_tf(x)=f(x)$  for each $x\in F$ if
$f$ is continuous and vanishes at infinity.
Our main theorem can be stated in terms of processes as follows.

\begin{corollary}\label{C1.3}
If $X$ is a continuous non-degenerate symmetric strong Markov process
which is a Feller process,  whose state space is $F$,  and whose Dirichlet
form  is invariant with respect to the
local symmetries of $F$, then the law of $X$ under $\bP^x$ is uniquely
defined, up to scalar multiples of the time parameter, for each 
$x\in F$.  
\end{corollary}

\begin{remark}\label{rem:osada}{\rm 
Osada \cite{Osada2006} 
constructed  diffusion processes on GSCs which are different from the 
ones considered here. While his processes are invariant 
with respect to some of the local isometries of the GSC, they are
not invariant with respect to the full set of local isometries. }
\end{remark}

In Section 2 we give precise definitions, introduce the notation we
use, and prove some preliminary lemmas. In Section 3 we prove 
Proposition \ref{ebbkz}. In Section 4 we develop the properties
of Dirichlet forms $\sE \in \fE$, and 
in Section 5 we prove Theorem~\ref{tmain}.

The idea of our proof is the following. The main work is showing that
if $\sA, \sB$ are any two Dirichlet forms
in $\fE$, then they are comparable. 
(This means that $\sA$ and $\sB$ have the same domain $\sF$,
and that there exists a constant $c=c(\sA,\sB)>0$ such that
$c \sA(f,f) \le \sB(f,f) \le c^{-1} \sA(f,f)$ for $f \in \sF$.)
We then let $\lam$ be the
largest positive real such that $\sC=\sA-\lam \sB\geq 0$. If $\sC$
were also in $\fE$, then $\sC$ would be comparable to $\sB$, and so
there would exist $\varepsilon>0$ such that $\sC-\varepsilon
\sB\geq 0$, contradicting the definition of $\lam$. In fact we
cannot be sure that $\sC$ is closed, so instead we consider
$\sC_\delta=(1+\delta)\sA-\lam \sB$, which is easily seen to be in
$\fE$. We then need uniform estimates in $\delta$ to obtain a
contradiction.

To show $\sA, \sB\in \fE$ are comparable requires heat kernel estimates
for an arbitrary element of $ \fE$. 
Using symmetry arguments as in \cite{BB4},
we show that the  estimates for corner moves and slides and the 
coupling argument of \cite[Section 3]{BB4}
can be modified so as to apply to any element $\sE \in \fE$. 
It follows that the elliptic Harnack inequality holds for any such $\sE$. 
Resistance arguments, as in \cite{BB3, McG},
combined with results in \cite{GT3} then lead to the desired
heat kernel bounds.
(Note that the results of \cite{GT3} that we use are also available in \cite{BBKT}.) 

A key point here is that the constants in the Harnack inequality, and
consequently also the heat kernel bounds, only depend on the GSC $F$,
and not on the particular element of $\fE$.  This means that we need
to be careful about the dependencies of the constants.

The symmetry arguments are harder than in 
\cite[Section 3]{BB4}. In \cite{BB4} the approximating
processes were time changed reflecting Brownian motions, and the 
proofs used the convenient fact that a
reflecting Brownian motion in a Lipschitz domain in $\bR^d$
does not hit sets of dimension $d-2$. Since we do not have such 
approximations for the processes corresponding to an arbitrary element 
$\sE \in \fE$, we have to work with the diffusion $X$ associated with $\sE$,
and this process might hit sets of
dimension $d-2$. 
(See \cite[Section 9]{BB4} for examples of GSCs in dimension $3$ for which
the process $X$ hits not just lines but also points.)

We use $C_i$ to denote finite positive constants which depend only on the
GSC, but which may change between each appearance.  Other finite positive
constants will be written as $c_i$.

\section{Preliminaries}\label{SP}

\subsection{Some general properties of  Dirichlet forms}

\ms We begin with a general result on local Dirichlet forms.
For definitions of local and other terms related to Dirichlet
forms, see \cite{FOT}.  Let $F$ be a compact metric space and $m$ a
Radon (i.e. finite) measure on $F$.  For any Dirichlet form $(\sE,\sF)$ on
$L^2(F,m)$
we define
\begin{equation}\label{DE1}
\sE_1(u,u)=\sE(u,u)+\|u\|_2^2. 
\end{equation}
Functions in $\sF$ are only defined up to quasi-everywhere equivalence
(see \cite{FOT} p. 67); we use a quasi-continuous modification of
elements of $\sF$ throughout the paper.
We write $\langle \cdot, \cdot \rangle$ for the inner product in
$L^2(F,m)$ and $\langle \cdot, \cdot \rangle_S$ for the inner product
in a subset $S\subset F$. 

\begin{theorem}\label{dfgen}
Suppose that $(\sA,\sF)$, $(\sB,\sF)$ are 
local regular conservative
irreducible Dirichlet forms on $L^2(F,m)$ and that 
\begin{equation} \label{e:aleb}
\sA(u,u) \le \sB(u,u) \q \hbox{ for all } u \in \sF. 
\end{equation}
Let $\delta>0$, and  $\sE=(1+\delta) \sB -\sA$. Then $(\sE,\sF)$ is a regular
local conservative irreducible Dirichlet form on $L^2(F,m)$.
\end{theorem}

\proof
It is clear that $\sE$ is a non-negative symmetric form,
and is local.

To show that $\sE$ is closed, let $\{u_n\}$ be a Cauchy sequence with
respect to $\sE_1$.  Since $\sE_1(f,f) \ge (\delta\wedge1) \sB_1(f,f)$,
$\{u_n\}$ is a Cauchy sequence with respect to $\sB_1$. Since $\sB$ is
a Dirichlet form and so closed, there exists $u \in \sF$ such that
$\sB_1(u_n-u,u_n-u) \to 0$. As $\sA \le \sB$ we have $\sA(u_n-u,u_n-u)
\to 0$ also, and so $\sE_1(u_n-u,u_n-u) \to 0$, proving that
$(\sE,\sF)$ is closed.

Since $\sA$ and $\sB$ are conservative and $F$ is compact,
$1\in \sF$ and $\sE(1,h)=0$ for all $h\in \sF$,
which shows that $\sE$ is conservative by
\cite[Theorem 1.6.3 and Lemma 1.6.5]{FOT}.

We now show that $\sE$ is Markov.
By \cite[Theorem 1.4.1]{FOT}  it is enough to prove that
$\sE({\bar u},{\bar u}) \le \sE(u,u)$ for $u \in \sF$, where
we let ${\bar u}=0\vee (u\wedge 1)$. 
Since $\sA$ is local and $u_+ u_-=0$, we have $\sA(u_+,u_-)=0$
(\cite[Proposition 1.4]{Schum}). Similarly $\sB(u_+,u_-)=0$, giving 
$\sE(u_+,u_-)=0$. Using this, we have 
\begin{equation}\label{Eplusminus}
\sE(u,u)=\sE(u_+,u_+)-2\sE(u_+,u_-)+\sE(u_-,u_-)\ge \sE(u_+,u_+) 
\end{equation} 
for $u \in \sF$. 
Now let $v=1-u$. Then ${\bar u}=(1-v_+)_+$ , so
\begin{align*}
\sE(u,u)&=\sE(v,v)
\ge \sE(v_+,v_+)=\sE(1-v_+,1-v_+)\\
&\ge \sE((1-v_+)_+,(1-v_+)_+)=\sE({\bar u},{\bar u}),
\end{align*}
and hence $\sE$ is Markov.

As $\sB$ is regular, it has a core $\sC \subset \sF$.
Let $u \in \sF$. As $\sC$ is a core for $\sB$,
there exist $u_n \in \sC$ such that
$\sB_1(u-u_n,u-u_n) \to 0$. Since  $\sA \le \sB$,
 $\sA_1(u_n-u,u_n-u) \to 0$ also, and so
 $\sE_1(u_n-u,u_n-u) \to 0$. Thus $\sC$ is
dense in $\sF$ in the $\sE_1$ norm  
(and it is dense in $C(F)$ in the supremum norm since it is a 
core for $\sB$), so $\sE$ is regular.

Let $A\subset F$ be invariant for the semigroup
corresponding to $\sE$. By \cite[Theorem 1.6.1]{FOT}, this is equivalent to
the following: $1_Au\in \sF$ for all $u\in \sF$ and
\begin{equation}\label{invfot}
\sE(u,v)=\sE(1_Au,1_Av)+\sE(1_{F-A}u,1_{F-A}v)\qquad\forall u,v\in \sF.
\end{equation}
Once we have $1_Au\in \sF$, since $(1_A u)(1_{F-A} u)=0$ we have
$\sA( 1_A u, 1_{F-A} u)=0$, and we obtain 
(\ref{invfot}) for $\sA$ also. 
Using \cite[Theorem 1.6.1]{FOT} again, we see that $A$
is invariant for the semigroup corresponding to $\sA$. Since $\sA$ is
irreducible, we conclude that either $m(A)=0$ or $m(X-A)=0$ holds and
hence that $(\sE,\sF)$ is irreducible. \qed

\begin{remark} 
{\rm This should be compared with the situation for
Dirichlet forms on finite sets, which is the context of the
uniqueness results in \cite{Me1, Sab1}. In that case the
Dirichlet forms are not local, and given $\sA$, $\sB$ satisfying
\eqref{e:aleb} there may exist $\delta_0>0$ such that
$(1+\delta)\sB -\sA$ fails to be a Dirichlet form for
$\delta\in (0,\delta_0)$. }
\end{remark}

For the remainder of this section we assume that $(\sE,\sF)$ is a
local regular Dirichlet form on $L^2(F,m)$, that $1\in \sF$ and
$\sE(1,1)=0$.  We write $T_t$ for the semigroup associated with $\sE$,
and $X$ for the associated diffusion.

\begin{lemma}\label{Lemma2.2}  $T_t$ is recurrent and conservative.
\end{lemma}

\proof $T_t$  is recurrent by \cite[Theorem 1.6.3]{FOT}.
Hence by \cite[Lemma 1.6.5]{FOT}  $T_t$ is
conservative. \qed

\sms Let $D$ be a Borel subset of $ F$.  We write $T_D$ for the hitting time of $D$, and
$\tau_D$ for the exit time of $D$:
\begin{equation}\label{TD}
 T_D= T_D^X=\inf\{ t \ge 0: X_t \in D\}, \qq 
\tau_D= \tau_D^X = \inf\{t \ge 0: X_t \not\in D\}. 
\end{equation}
Let $\ol T_t$ be the semigroup of $X$ killed on exiting $D$, and
$\ol X$ be the killed process. Set
$$ q(x) = \bP^x( \tau_D = \infty), $$
and
\begin{equation}\label{defED} 
E_D =\{x : q(x)=0 \}, \qq Z_D=\{x : q(x)=1 \}. 
\end{equation}

\begin{lemma}\label{Lemma2.3} 
Let $D$ be a Borel subset of $F$.
Then $\ m(D - (E_D \cup Z_D))=0$. Further, $E_D$ and $Z_D$ are
invariant sets for the killed process $\ol X$,
and $Z_D$ is invariant for $X$.
\end{lemma}

\proof If $f \ge 0$,
$$ \langle \ol T_t (f 1_{E_D}), 1_{D-E_D} q \rangle 
 =   \langle f 1_{E_D}, \ol T_t (1_{D-E_D} q) \rangle 
\le  \langle f 1_{E_D}, \ol T_t q \rangle =0. $$
So  $\ol T_t (f 1_{E_D})=0$ on $D-E_D$
and hence (see \cite[Lemma 1.6.1(ii)]{FOT}) 
$E_D$ is invariant for $\ol X$. 

Let $A=\{x:P^x(\tau_D<\infty)>0\}=Z_D^c$.  The set 
$A$ is an invariant set of the process $X$ by \cite[Lemma 4.6.4]{FOT}. 
Using the fact that $\ol X = X$, $\bP^x$-a.s. for 
$x \in Z_D$ and  \cite[Lemma 1.6.1(ii)]{FOT}, we see that 
$A$ is an invariant set of the process $\ol X$ as well. 
So we see that $Z_D$ is invariant both for $X$ and $\ol X$. 
In order to 
prove $m(D - (E_D \cup Z_D))=0$, it suffices to show that
$\bE^x[\tau_D]<\infty$ for a.e. $x\in A\cap D$. 
Let $U_D$ be
the resolvent of the killed process $\ol X$.   Since $A\cap D$ is of finite
measure, the proof of Lemma 1.6.5 or Lemma 1.6.6 of \cite{FOT} give 
$U_D1(x)<\infty$ for a.e. $x\in A\cap D$, so we obtain $\bE^x[\tau_D]<\infty$. 
\qed

Note that in the above proof 
we do not use the boundedness of $D$, but only the fact that $m(D)<\infty$.

\medskip
Next, we give some general facts on harmonic and caloric functions. 
Let $D$ be a Borel subset in $F$ and let $h : F \to \bR$. 
There are two possible definitions of $h$ being harmonic
in $D$.
The probabilistic one is that $h$ is harmonic in $D$ 
if $h(X_{t\land \tau_{D'}})$ is a 
uniformly integrable martingale
under $\bP^x$ for q.e.\ $x$ whenever $D'$ is a relatively open subset of $D$. 
The Dirichlet form definition is
that $h$ is harmonic with respect to $\sE$ in $D$ 
if $h \in \sF$
and $\sE(h,u)=0$ whenever $u\in \sF$ is continuous
and the support of $u$ is contained in $D$.

The following is well known to experts. 
We will use it in the proofs 
of Lemma~\ref{equal-hit} and  Lemma~\ref{enorm}. 
(See \cite{Chnew} 
for the  equivalence of the two notions of harmonicity in a very 
general framework.)
Recall that
$\bP^x( \tau_D < \infty)=1$ for  $x\in E_D$.

 \begin{proposition}\label{s2p2} \, (a) \, 
Let $(\sE,\sF)$ and $D$ satisfy the above conditions, and let $h\in \sF$ be bounded. 
Then $h$ is harmonic in a domain $D$ in the probabilistic
sense if and only if it is harmonic in the Dirichlet form sense. \\
(b) \, If $h$ is a bounded Borel measurable function  in $D$ and $D'$ is a relatively open 
subset of $D$, then $h(X_{t\land \tau_{D'}})$ is a martingale
under $\bP^x$ for q.e.\ $x\in E_D$ 
if and only if 
$h(x)=\bE^x[h(X_{\tau_{D'}})]$ for q.e.\ $x\in E_D$. 
\end{proposition}

\proof (a) By \cite[Theorem 5.2.2]{FOT}, we have the 
Fukushima decomposition 
$h(X_{t})-h(X_0)=M^{[h]}_t+N^{[h]}_t,$
where $M^{[h]}$ is a square integrable martingale additive functional 
of finite energy and $N^{[h]}$ is a continuous additive functional having zero energy
(see \cite[Section 5.2]{FOT}). 
We need to  consider the Dirichlet form $(\sE,\sF_D)$ 
where $\sF_{D}= \{ f \in \sF: \supp(f) \subset D \}$,  
and denote the corresponding semigroup as $P_t^D$.

If $h$ is harmonic in the Dirichlet form sense, then by 
the discussion in \cite[p.~218]{FOT} and
\cite[Theorem 5.4.1]{FOT}, 
we have 
$\bP^x(N^{[h]}_t=0,~\forall t< \tau_D)=1$ q.e. $x\in F$. Thus, 
$h$ is harmonic in the probabilistic sense. Here the notion of the spectrum from \cite[Sect.~2.3]{FOT} and especially \cite[Theorem~2.3.3]{FOT} are used. 

To show that  being harmonic in the probabilistic sense  implies being 
harmonic in the Dirichlet form sense is the  delicate part of this proposition. 
Since $Z_D$ is $P^D_t$-invariant (by Lemma~\ref{Lemma2.3}) and 
$h(X_t)$ is a bounded martingale under $\bP^x$ for $x\in Z_D$, we have
\[P^D_t(h1_{Z_D})(x)=1_{Z_D}(x)P^D_th(x)=1_{Z_D}(x)E^x[h(X_t)]=h1_{Z_D}(x).\]
Thus by \cite[Lemma 1.3.4]{FOT}, we have $h1_{Z_D}\in \sF$ and 
$\sE(h1_{Z_D},v)=0$ for all $v\in \sF$. 
Next, note that on $Z_D^c$ we have $H_Bh=h$, according to the definition of $H_B$ on 
page 150 of \cite{FOT} and Lemma~\ref{Lemma2.3}, which implies 
$H_B\left(h1_{Z_D^c}\right)=h1_{Z_D^c}$. Then from \cite[Theorem~4.6.5]{FOT}, 
applied with $\widetilde u=h1_{Z_D^c}=h-h1_{Z_D}\in \sF$ 
and $B^c=D$, we conclude that $h1_{Z_D^c}$ 
is harmonic in the Dirichlet form sense. 
Thus $h=h1_{Z_D^c}+h1_{Z_D}$ is harmonic in the Dirichlet form sense in $D$.

\ms
\noindent(b) If $h(X_{t\land \tau_{D'}})$ is a martingale
under $\bP^x$ for q.e.\ $x\in E_D$, then 
$\bE^x[h(X_{s\land \tau_{D'}})]=\bE^x[h(X_{t\land \tau_{D'}})]$ 
for q.e. $x\in E_D$ and for all
$s,t\ge 0$, where we can take $s\downarrow 0$ and $t\uparrow \infty$ and interchange the
limit and the expectation since 
$h$ is bounded. Conversely, if $h(x)=E^x[h(X_{\tau_{D'}})]$ 
for q.e.\ $x\in E_D$, then by the strong Markov property,
$h(X_{t\land \tau_{D'}}) 
=E^x[h(X_{\tau_{D'}})|{\cal F}_{t\wedge \tau_{D'}}]$ under $\bP^x$ 
for q.e. $x\in E_D$, so $h(X_{t\land \tau_{D'}})$ is a martingale 
under $\bP^x$ for q.e.\ $x\in E_D$. 
\qed

We call a function $u:\bR_+\times F\to\bR$ caloric in $D$ in the
probabilistic sense if $u(t,x)=\bE^x[f(X_{t\land \tau_{D}})]$ for some
bounded Borel $f:F\to\bR$. It is natural to view $u(t,x)$ as the
solution to the heat equation with boundary data defined by $f(x)$
outside of $D$ and the initial data defined by $f(x)$ inside of
$D$. We call a function $u:\bR_+\times F\to\bR$ caloric in $D$ in the
Dirichlet form sense if there is a function $h$ which is harmonic in
$D$ and a bounded Borel $f_D:F\to\bR$ which vanishes outside of $D$
such that $u(t,x)=h(x)+\ol T_t f_D$. Note that $\ol T_t$ is the
semigroup of $X$ killed on exiting $D$, which can be either defined
probabilistically as above or, equivalently, in the Dirichlet form
sense by Theorems 4.4.3 and A.2.10 in~\cite{FOT}.
  
\begin{proposition}\label{p-caloric} 
Let $(\sE,\sF)$ and $D$ satisfy the above conditions, and let $f\in
\sF$ be bounded and $t\ge0$.  Then $$\bE^x[f(X_{t\land
\tau_{D}})]=h(x)+\ol T_t f_D$$ q.e., where $h(x)=\bE^x[f(X_{
\tau_{D}})]$ is the harmonic function that coincides with $f$ on
$D^c$, and $f_{D}(x)=f(x)-h(x)$.
\end{proposition}

\proof By Proposition~\ref{s2p2}, $h$ is uniquely defined in
the probabilistic and Dirichlet form senses, and
$h(x)=\bE^x[h(X_{t\land \tau_{D}})]$. Note that $f_{D}(x)$ vanishes
q.e. outside of $D$. Then we have $\bE^x[f_D(X_{t\land \tau_{D}})]=\ol
T_t f_D$ by Theorems 4.4.3 and A.2.10 in \cite{FOT}.  \qed

Note that the condition $f\in \sF$ can be relaxed (see the proof of
Lemma~\ref{equal-hit}).

\sms
We show a general property of local Dirichlet
forms which will be used in the proof of Proposition~\ref{P2}.  Note
that it is not assumed that $\sE$ admits a {\it carr\'{e} du champ}.  
Since $\sE$ is regular, $\sE(f,f)$ can be written in terms 
of a measure $\Gam({f,f})$, 
the energy measure of $f$, as follows.  Let $\sF_b$ be the elements
of $\sF$ that are essentially bounded. If $f\in \sF_b$, then
$\Gam({f,f})$ is defined to be the unique
smooth Borel measure on $F$ satisfying
$$\int_F{g}d\Gam({f,f})=2\sE (f,fg)-\sE (f^2,g),\qquad  g\in
\sF_b.$$

\begin{lemma}\label{lem-zem}
If $\sE$ is a local regular Dirichlet form with domain $\sF$, then for any 
$f\in\sF\cap L^\infty(F)$ 
we have $\Gam(f,f)(A)=0$, where $A=\{x\in F:f(x)=0\}$.  
\end{lemma}

\proof Let $\sigma^f$ be the measure on  $\bR$ which is the image 
 of the measure $\Gam(f,f)$  on $F$ under the function $f:F\to\bR$. 
By \cite[Theorem 5.2.1, Theorem 5.2.3]{BH} and the chain rule, 
$\sigma^f$ is absolutely continuous with respect to one-dimensional 
Lebesgue measure on $\bR$.  Hence  $\Gam(f,f)(A)=\sigma^f(\{0\})=0$. \qed

\begin{lemma}\label{lem-Feller}
Given a $m$-symmetric Feller process on $F$, the corresponding Dirichlet form 
$(\sE,\sF)$ is  regular.  
\end{lemma}

\proof 
First, we note the following: if $H$ is dense in $L^2(F,m)$, then $U^1 (H)$ 
is dense in $\sF$, where $U^1$ is the $1$-resolvent operator. 
This is because
$U^1: L^2\to \sD(\sL)$ is an isometry where the norm of $g\in\sD(\sL)$ 
is given by $\|g\|_{\sD(\sL)}:=\|(I-\sL)g\|_2$, and $\sD(\sL)
\subset \sF$ is a continuous dense embedding (see, for example
\cite[Lemma 1.3.3(iii)]{FOT}). 
Here $\sL$ is the generator corresponding to $\sE$. 
Since $C(F)$ is dense in $L^2$ and $U^1 (C(F))\subset \sF\cap C(F)$
as the process is Feller, we see that $\sF\cap C(F)$ is dense in 
$\sF$ in the $\sE_1$-norm.

Next we need to show that $u\in C(F)$ can be approximated 
with respect to the supremum norm by functions in $\sF\cap C(F)$. This is easy,
since $T_tu\in \sF$ for each $t$, is continuous since we have a Feller
process, and $T_tu\to u$ uniformly by \cite[Lemma III.6.7]{RW}.
\qed

\begin{remark}{\rm The proof above uses the fact that $F$ is
compact.
However, it can be easily generalized to a Feller process on a 
locally compact separable metric space by a
standard truncation argument -- for example by using 
\cite[Lemma 1.4.2(i)]{FOT}.  
}\end{remark}

\subsection{Generalized Sierpinski carpets}\label{subGSC}

Let $d\geq 2$, $F_0=[0,1]^d$, and let 
$\LF \in \bN$, $\LF\geq 3$, be fixed. For
$n \in \bZ$ let $\sQ_n$ be the collection of closed cubes of side
$\LF^{-n}$ with vertices in $\LF^{-n}\bZ^d$. For $A\subseteq\bR^d$,
set
$$\sQ_n(A)=\{Q \in \sQ_n: \inter(Q) \cap A \neq \emptyset \} .$$
For $Q \in\sQ_n$, let $\Psi_Q$ be the orientation preserving affine
map  (i.e.\ similitude with no rotation part) which maps $F_0$ onto
$Q$. We now define a decreasing sequence $(F_n)$ of closed subsets
of $F_0$. Let $1\leq \MF \le \LF^d$ be an integer, and let $F_1$ be
the union of $\MF$ distinct elements of $\sQ_1(F_0)$. We impose the
following conditions on $F_1$. 

\begin{description}
\item{(H1)} (Symmetry) $F_1$ is preserved by all the isometries
of the unit cube $F_0$.

\item{(H2)} (Connectedness) ${\rm Int}(F_1)$ is connected.

\item{(H3)} (Non-diagonality) Let $m \ge 1$
and $B\subset F_0$ be a cube of side length
$2\LF^{-m}$, which is the union of $2^d$ distinct elements of
$\sQ_m$.  
Then if $\inter(F_1 \cap B)$ is non-empty, it is connected.

\item{(H4)} (Borders included) $F_1$ contains the line segment $\{x: 0\leq x_1
\leq 1, x_2= \cdots =x_d=0 \}$.
\end{description}

\medskip
We may think of $F_1$ as being derived from $F_0$ by removing the
interiors of $\LF^d - \MF$ cubes in $\sQ_1(F_0)$. Given $F_1$, $F_2$
is obtained by removing the same pattern from each of the cubes in
$\sQ_1(F_1)$. Iterating, we obtain a sequence $\{F_n\}$, where $F_n$
is the union of $\MF^n$ cubes in $\sQ_n(F_0)$. Formally, we define
$$F_{n+1}=\bigcup_{Q\in\sQ_n(F_n)}\Psi_Q(F_1)=
\bigcup_{Q\in\sQ_1(F_1)}\Psi_Q(F_n), \q n\geq 1.$$
 We call the set $F= \cap_{n=0}^\infty F_n$ a  generalized Sierpinski
  carpet (GSC).
The Hausdorff dimension of $F$ is $d_f=d_f(F)=\log \MF/\log \LF$.
Later on we will also discuss the unbounded GSC
$\wt F=\cup_{k=0}^\infty \LF^k F$,
where $rA=\{rx:x\in A\}$.

\medskip
Let
$$\mu_n (dx) = (\LF^d/\MF)^{n} 1_{F_n} (x) dx,$$
and let $\mu$ be the weak limit of the $\mu_n$; $\mu$ is a constant
multiple of the Hausdorff $x^{d_f}$ - measure on $F$.
For $x, y \in F$ we write $d(x,y)$ for the length of the shortest path
in $F$ connecting $x$ and $y$. 
Using (H1)--(H4) we have that $d(x,y)$ is
comparable with the Euclidean distance $|x-y|$. 

\sms \begin{remark}\label{rem:correcterror} {\rm 

1. There is an error in \cite{BB4}, where it was only
assumed that (H3) above holds when $m=1$. However, that assumption 
is not strong enough to imply the connectedness of the set $J_k$ 
in \cite[Theorem 3.19]{BB4}.
To correct this error, we replace the (H3) in \cite{BB4} by the (H3) 
in the current paper. \\
2. The {\sl standard SC} in dimension $d$ is the GSC with $L_F=3$,
$\MF=3^d-1$, and with $F_1$ obtained from $F_0$ by removing the middle
cube. We have allowed $\MF=L_F^d$, so that our GSCs do include
the `trivial' case  $F=[0,1]^d$. 
The `Menger sponge' (see the picture on \cite{Man}, p. 145) is
one example of a GSC, and has $d=3$, $L_F=3$, $\MF=20$.  
}\end{remark}

\begin{figure}
\begin{center}
\begin{picture}(120,120)(0,0)\small\thinlines
\setlength{\unitlength}{40pt}
\put(0,0){\line(0,1){3}\line(1,0){3}}
\put(3,3){\line(0,-1){3}\line(-1,0){3}}
\put(1,1){\line(0,1){1}\line(1,0){1}}
\put(2,2){\line(0,-1){1}\line(-1,0){1}}
\linethickness{.25pt}
\multiput(0,0)(.025,0){40}{\line(0,1){1}}
\linethickness{2.5pt}
\qbezier[10](1,0)(1,0.5)(1,1)
\put(0,1){\linethickness{.25pt}
\multiput(0,0)(.025,0){40}{\line(0,1){1}}
\linethickness{2.5pt}
\qbezier[10](0,1)(0.5,1)(1,1)}
\put(1,1){\rule{40pt}{40pt}}
\put(0,0){\put(.3333,.3333){\rule{13.3333pt}{13.3333pt}}}
\put(0,1){\put(.3333,.3333){\rule{13.3333pt}{13.3333pt}}}
\put(0,2){\put(.3333,.3333){\rule{13.3333pt}{13.3333pt}}}
\put(1,0){\put(.3333,.3333){\rule{13.3333pt}{13.3333pt}}}
\put(2,0){\put(.3333,.3333){\rule{13.3333pt}{13.3333pt}}}
\put(2,1){\put(.3333,.3333){\rule{13.3333pt}{13.3333pt}}}
\put(2,2){\put(.3333,.3333){\rule{13.3333pt}{13.3333pt}}}
\put(1,2){\put(.3333,.3333){\rule{13.3333pt}{13.3333pt}}}
\end{picture}
\hfil 
\begin{picture}(120,120)(0,0)\small\thinlines
\setlength{\unitlength}{40pt}
\put(0,0){\line(0,1){3}\line(1,0){3}}
\put(3,3){\line(0,-1){3}\line(-1,0){3}}
\put(1,1){\line(0,1){1}\line(1,0){1}}
\put(2,2){\line(0,-1){1}\line(-1,0){1}}
\linethickness{.25pt}
\multiput(0,0)(.025,0){40}{\line(0,1){1}}
\linethickness{2.5pt}
\qbezier[10](1,0)(1,0.5)(1,1)
\qbezier[10](0,0)(0,0.5)(0,1)
\qbezier[10](0,0)(0.5,0)(1,0)
\put(0,1){\linethickness{.25pt}
\multiput(0,0)(.025,0){40}{\line(0,1){1}}
\linethickness{2.5pt}
\qbezier[10](1,0)(1,0.5)(1,1)
\qbezier[10](0,1)(0.5,1)(1,1)
\qbezier[10](0,0)(0,0.5)(0,1)}
\put(1,1){\rule{40pt}{40pt}}
\put(0,0){\put(.3333,.3333){\rule{13.3333pt}{13.3333pt}}}
\put(0,1){\put(.3333,.3333){\rule{13.3333pt}{13.3333pt}}}
\put(0,2){\put(.3333,.3333){\rule{13.3333pt}{13.3333pt}}}
\put(1,0){\put(.3333,.3333){\rule{13.3333pt}{13.3333pt}}}
\put(2,0){\put(.3333,.3333){\rule{13.3333pt}{13.3333pt}}}
\put(2,1){\put(.3333,.3333){\rule{13.3333pt}{13.3333pt}}}
\put(2,2){\put(.3333,.3333){\rule{13.3333pt}{13.3333pt}}}
\put(1,2){\put(.3333,.3333){\rule{13.3333pt}{13.3333pt}}}
\end{picture}
\end{center}
\caption{Illustration for  Definition~\ref{dSn} in the 
case of the standard Sierpinski carpet and $n=1$. 
Let $A$ be the shaded set. The thick dotted lines 
give $\inter_F A$ on the left, and $\inter_r A$ on the right.
}\label{fig1}\end{figure}
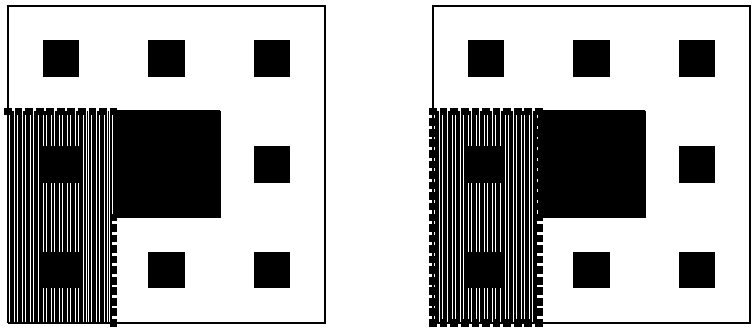

\begin{definition}{\rm\label{dSn} 

\sms Define:
$$ \sS_n=\sS_n(F) =\{ Q\cap F: Q \in \sQ_n(F)\} . $$

We will need to consider two different types of interior and boundary
for subsets of $F$ which consist of unions of  elements of $\sS_n$. First, for
any $A \subset F$ we write $\inter_F(A)$ for the  interior
of $A$ with respect to  the metric space $(F,d)$, 
and $\pd_F(A) = \ol A - \inter_F(A)$.
Given any $U \subset \bR^d$ we write $U^o$ for the interior of $U$ in
with respect to the usual topology on $\bR^d$, and $\pd U=\ol U - U^o$
for the usual boundary of $U$.
Let $A$ be a finite union of  elements of $\sS_n$, so 
that $A = \cup_{i=1}^k S_i$,
where $S_i = F \cap Q_i$ and $Q_i \in \sQ_n(F)$. Then we define
$\inter_r (A) = F \cap((\cup_{i=1}^k Q_i)^o)$, and $\pd_r(A)= A - \inter_r (A)$.
We have $\inter_r(A) = A - \pd( \cup_{i=1}^k Q_i)$.
(See Figure~\ref{fig1}).

}\end{definition}

\begin{definition}{\rm\label{dphi}
We define the folding map $\vp_S: F\to S$ for $S \in \sS_n(F)$ as follows.
Let $\ol \vp_0: [-1,1]\to \bR$ be defined by 
$\ol \vp_0(x)=|x|$ for  $|x|\leq 1$,
and then extend the domain of $\ol \vp_0$ to all of $\bR$ by periodicity, so
that $\ol \phi_0(x+2 n)= \ol \phi_0(x)$ for all $x \in \bR$, $n \in \bZ$.
If $y$ is the point of $S$ closest to the origin, define
$\vp_S(x)$ for $x\in F$ to be the point whose $i^{th}$
coordinate is $y_i+ \LF^{-n} \ol\vp_0(\LF^n(x_i-y_i))$.

\sms It is straightforward to check the following

\begin{lemma}\label{vp-prop}
(a) $\vp_S$ is the identity on $S$ and for each $S' \in \sS_n$,
$\vp_S : S' \to S$ is an isometry. \\
(b) If $S_1, S_2 \in \sS_n$ then
\begin{equation} \label{phi-p1}
 \phi_{S_1} \circ \phi_{S_2} = \phi_{S_1}. 
\end{equation}
(c) Let $x,y \in F$. If there exists $S_1 \in \sS_n$ such that 
$\vp_{S_1}(x) = \vp_{S_1}(y)$, then $\vp_{S}(x) = \vp_{S}(y)$ for
every $S \in \sS_n$.\\
(d) Let $S \in \sS_n$ and $S'\in \sS_{n+1}$. If $x, y \in F$
and  $\vp_{S}(x) = \vp_{S}(y)$ then $\vp_{S'}(x) = \vp_{S'}(y)$.
\end{lemma}

\sms
Given $S \in \sS_n$, $f: S \to \bR$ and $g: F \to \bR$ 
we define the unfolding and restriction operators by
$$ U_S f = f \circ \phi_S, \q R_S g = g|_S. $$
Using \eqref{phi-p1}, we have that if $S_1, S_2 \in \sS_n$ then 
\begin{equation}\label{URproj}
 U_{S_2} R_{S_2} U_{S_1} R_{S_1} =  U_{S_1} R_{S_1}. 
\end{equation}}\end{definition}

\begin{definition}\label{scalefac}
{\rm  
We define the {\em length} and {\em mass} scale factors of
$F$ to be $L_F$ and $\MF$ respectively. 

Let $D_n$ be the network of diagonal crosswires
obtained by joining each vertex of a cube $Q \in \sQ_n$ to the
vertex at the center of the cube by a wire of unit resistance
-- see \cite{BB3, McG}. Write $R_n^D$ for the resistance 
across two opposite faces of $D_n$. Then it is proved
in  \cite{BB3, McG} that there exists $\rho_F$ such 
that there exist constants $C_i$, depending only on the
dimension $d$, such that
\begin{equation}\label{rhodef}
 C_1 \rho_F^n \le R_n^D  \le C_2 \rho_F^n .
\end{equation}
We remark that $\rho_F \le L_F^2/\MF$ -- 
see \cite[Proposition 5.1]{BB4}.
}\end{definition}

\medskip

\subsection{$F$-invariant Dirichlet forms}\label{subFinv}

Let $(\sE, \sF)$ be a local regular Dirichlet form on $L^2(F,\mu)$.
Let $S \in \sS_n$. We set
\begin{equation}\label{defES}
\sE^{S}(g,g)=\frac{1}{\MF^n} \sE(U_{S} g, U_{S} g).
\end{equation}
and define the domain of $\sE^{S}$ to be 
$\sF^{S} = \{ g: g\mbox{ maps }S \mbox{ to } \bR, U_S g \in \sF \}$.
We write $\mu_S = \mu |_S$.

\begin{definition} \label{deffE} {\rm
Let $(\sE,\sF)$ be a  Dirichlet form on $L^2(F,\mu)$.
We say that $\sE$ is an \emph{$F$-invariant Dirichlet form}
or that 
\emph{$\sE$ is invariant with respect to all the local symmetries of $F$}
if the following items (1)--(3) hold:
\begin{description}
\item{(1)} If $S \in \sS_n(F)$, then 
$U_{S}R_S f \in \sF$ (i.e. $R_S f \in \sF^S$) for any $f\in \sF$. 
\item{(2)}  Let $n \ge 0$ and $S_1, S_2$ be any two elements of $\sS_n$, and
let $\Phi$ be any isometry of $\bR^d$ which maps $S_1$ onto $S_2$. 
(We allow $S_1=S_2$.) 
If $f \in \sF^{S_2}$, then $f \circ \Phi \in \sF^{S_1}$ and
\begin{equation}\label{isom}
\sE^{S_1}( f \circ \Phi, f \circ \Phi)=\sE^{S_2}(f,f).
\end{equation} 
\item{(3)} For all $f\in \sF$
\begin{equation}\label{sumDF}
\sE(f,f)=\sum_{S \in \sS_n(F)} \sE^{S}(R_S f, R_S f).
\end{equation}

\end{description}

\noindent We write $\fE$ for the set of $F$-invariant, non-zero, 
local, regular, conservative  Dirichlet forms. }
\end{definition}

\begin{remark}\label{rem-2.3ss}{\rm  We cannot exclude at this point the 
possibility that the energy measure of $\sE \in \fE$ may charge the
boundaries of cubes in $\sS_n$. See Remark \ref{rem:chargeboundary}.
} \end{remark}

We will not need the following definition of scale invariance
until we come to the proof of Corollary \ref{C1.2} in Section
\ref{Uni}.  

\begin{definition}\label{def-scaleinv} {\rm
Recall that $\Psi_Q$, $Q \in \sQ_1(F_1)$ are the similitudes which 
define $F_1$. Let $(\sE,\sF)$ be a Dirichlet form on $L^2(F,\mu)$
and suppose that
\begin{equation}\label{psi-dom}
 f \circ \Psi_Q \in \sF \, \hbox{ for all } Q \in \sQ_1(F_1), \, f \in \sF.
\end{equation}
Then we can define the {\em replication} of $\sE$ by
\begin{equation}\label{def-rep}
 \sR \sE(f,f) = \sum_{Q \in \sQ_1(F_1)} \sE(f \circ \Psi_Q, f \circ \Psi_Q).
\end{equation}
We say that $(\sE,\sF)$ is {\em scale invariant} if \eqref{psi-dom}
holds, and there exists $\lam>0$ such that $\sR \sE = \lam \sE$.}
\end{definition}

\begin{remark} 
{\rm We do not have any direct proof that if $\sE \in \fE$
then \eqref{psi-dom} holds. Ultimately, however, this will follow
from Theorem \ref{tmain}. }
\end{remark}

\begin{lemma}\label{corL}
Let $(\sA,\sF_1)$, $(\sB, \sF_2) \in\fE$ with $\sF_1=\sF_2$ 
and $\sA\ge\sB$. Then
$\sC=(1+\delta)\sA- \sB\in\fE$ for any $\delta>0$.
\end{lemma}

\proof It is easy to see that Definition \ref{deffE} holds.
This and Theorem \ref{dfgen} proves the lemma.  \qed

\begin{proposition}\label{P4}
If $\sE \in \fE$ and $S \in \sS_n(F)$,
then $(\sE^{S}, \sF^{S})$ is a local regular
Dirichlet form on $L^2(S,\mu_S)$.
\end{proposition}

\proof (Local): If $u,v$ are in $\sF^{S}$ with compact
support and $v$ is constant in a neighborhood of the support of $u$,
then $U_S u,U_S v$ will be in $\sF$, and by the local property
of $\sE$, we have $\sE(U_S u,U_S v)=0$. Then by \eqref{defES} we have
$\sE^{S}(u,v)=0$.

\ms\noi (Markov):
Given that $\sE^{S}$ is local, we have the Markov property
by the same proof as that in Theorem \ref{dfgen}.

\ms\noi (Conservative):
Since $1\in \sF$, $\sE^{S}(1,1)=0$ by \eqref{defES}.

\ms\noi (Regular): 
If $h \in \sF$ then by \eqref{sumDF} 
$\sE^{S}( R_S h, R_S h) \le \sE(h,h)$. Let $f \in \sF^S$, so that
$U_S f \in \sF$. As $\sE$ is regular, given $\eps>0$ there exists
a continuous $g \in \sF$ such that 
$\sE_1(U_S f -g, U_S f -g) < \eps$.
Then $R_S U_S f - R_S g = f -R_S g$ on $S$, so
\begin{align*} 
\sE^{S}_1(f- R_S g,f-R_S g) &= \sE^{S}_1(R_S U_S f - R_S g, R_S U_S f - R_S g) \\
&\le \sE_1(U_S f -g, U_S f -g) < \eps.
\end{align*} 
As $R_S g$ is continuous, we see that $\sF^S\cap C(S)$ is dense in $\sF^S$ in the 
$\sE^S_1$ norm. One can similarly prove that $\sF^S\cap C(S)$ is dense in 
$C(S)$ in the supremum norm, so the regularity of $\sE^S$ is proved.

\ms\noi(Closed): 
If $f_m$ is Cauchy with respect to $\sE^{S}_1$, 
then $U_S f_m$ will be Cauchy with respect to $\sE_1$.
Hence $U_S f_m$ converges with respect to $\sE_1$, and it
follows that $R_S(U_S f_m)=f_m$ converges with respect to $\sE^{S}_1$.
\qed

Fix $n$ and define for functions $f$ on $F$
\begin{equation}\label{defPO}
\Theta f=\frac{1}{\MF^n} \sum_{S \in \sS_n(F)} U_S R_S f.
\end{equation}
Using \eqref{URproj} we have $\Theta^2=\Theta$, and so 
$\Theta$ is a projection operator. 
It is bounded on $C(F)$
and $L^2(F,\mu)$, and
moreover by \cite[Theorem 12.14]{Ru} is an orthogonal projection on
$L^2(F,\mu)$.  
Definition \ref{deffE}(1) implies that $\Theta: \sF \to \sF$.

\begin{proposition}\label{P2} Assume that $\sE$ is a local 
regular Dirichlet form on $F$, $T_t$ is its semigroup, and $U_S R_S
f\in \sF$ whenever $S\in\sS_n(F)$ and $f\in \sF$.  Then the following
are equivalent:
\begin{enumerate}
\item[{(a)}] \  For all $f\in \sF$, we have $\sE(f,f)=\sum_{S\in\sS_n(F)}
\sE^S(R_Sf,R_Sf)$;
\item[{(b)}] \ for all $f,g\in \sF$
\begin{equation}\label{sEQ}
\sE(\Theta f,g)=\sE(f,\Theta g);
\end{equation}
\item[{(c)}] \ $T_t\Theta f=\Theta T_tf$ a.e for any  $f\in L^2(F,\mu)$
and  $t\ge0$.
\end{enumerate}
\end{proposition}

\begin{remark} \rm Note that this proposition and the following
corollary do not use all the symmetries that are assumed in
Definition~\ref{deffE}(2). Although these symmetries are not needed
here, they will be essential later in the paper.
\end{remark}

\proof
To prove that $(a)\Rightarrow (b)$, note that 
(a) implies that
\begin{equation}\label{sEsum1}
 \sE(f,g) =  \sum_{T \in \sS_n(F)} \sE^T(R_T f, R_T g)=
\frac{1}{\MF^n} \sum_{T \in \sS_n(F)} \sE(U_T R_T f, U_T R_T g).
\end{equation}
Then using \eqref{defPO}, \eqref{sEsum1} and \eqref{URproj},
\begin{align*}
\sE(\Theta f,g) &= \frac{1}{\MF^n} \sum_{S \in \sS_n(F)} \sE( U_S R_S f,g) \\
&=  \frac{1}{\MF^{2n}} \sum_{S \in \sS_n(F)} \sum_{T \in \sS_n(F)} 
\sE( U_T R_T U_S R_S f, U_T R_T g) \\
&=  \frac{1}{\MF^{2n}} \sum_{S \in \sS_n(F)} \sum_{T \in \sS_n(F)} 
\sE(U_S R_S f, U_T R_T g).
\end{align*}
Essentially the same calculation shows that $\sE(f, \Theta g)$ is equal
to the last line of the above with the summations reversed.

Next we show that $(b)\Rightarrow (c)$.  
If $\sL$ is the generator corresponding to $\sE$,
$f\in\sD(\sL)$ and $g\in\sF$ then, writing $\langle f,g \rangle $ for 
$\int_F fg \, d\mu$, we have
$$\langle \Theta \sL f,g\rangle=\langle \sL f, 
\Theta g\rangle=-\sE(f,\Theta g)=-\sE(\Theta f,g)$$
by \eqref{sEQ} and the fact that $\Theta $ is self-adjoint in the $L^2$ sense.
By the definition of the generator corresponding to  a Dirichlet form, this is
equivalent to 
$$\Theta f\in\sD(\sL) \text{ \ and \ }\Theta \sL f=\sL \Theta f.$$ 
By \cite[Theorem 13.33]{Ru}, 
this implies that any bounded Borel function of $\sL$ 
commutes with $\Theta $. (Another good source on the
spectral theory of unbounded self-adjoint operators is \cite[Section VIII.5]{RS}.) 
In particular,  the $L^2$-semigroup $T_t$ of $\sL$ commutes with
 $\Theta $ in the $L^2$-sense.  This implies (c). 

In order to see that $(c)\Rightarrow (b)$, 
note that if $f,g\in \sF$, 
\begin{align*}
\sE(\Theta f,g)&=\lim_{t\to 0} t^{-1} \langle (I-T_t) \Theta f,g\rangle
=\lim t^{-1} \langle \Theta(I-T_t) f, g\rangle\\
&=\lim t^{-1} \langle (I-T_t)f,\Theta g\rangle=\lim t^{-1}
\langle f,(I-T_t) \Theta g\rangle\\
&=\sE(f, \Theta g).
\end{align*}

It remains to prove that $(b)\Rightarrow(a)$. This is the only
implication that uses the assumption that $\sE$ is local.  
It suffices to  assume $f$ and $g$ are bounded.

First, note the obvious relation
\begin{equation}\label{L213-E1}
\sum_{S\in\sS_n(F)}\frac{1_S(x)}{N_n(x)}=1
\end{equation}
 for any $x\in F$, where
\begin{equation}\label{Ndef}
N_n(x) = \sum_{S\in\sS_n(F)} 1_S(x)
\end{equation}
is the number of cubes $\sS_n$ whose interiors intersect $F$ and which
contain the point $x$.  We break the remainder of the proof into a
number of steps.

Step 1: We show that if $\Theta f=f$, then $\Theta (hf)=f(\Theta h)$.
To show this, we start with the relationship
$U_TR_TU_SR_Sf=U_SR_Sf$. Summing over $S\in \sS_n(F)$ and dividing
by $\MF^n$ yields
$$U_TR_Tf=U_TR_T\Theta (f)=\Theta f=f.$$
Since $R_S(f_1f_2)=R_S(f_1)R_S(f_2)$ and $U_S(g_1g_2)=U_S(g_1)U_S(g_2)$,
we have
$$\Theta(hf)=\frac{1}{\MF^n} \sum_{S\in \sS_n}(U_SR_Sf)(U_SR_Sh)
=\frac{1}{\MF^n} \sum_{S\in \sS_n} f(U_SR_Sh)=f(\Theta h).$$
In particular, $\Theta (f^2)=f\Theta f=f^2$.

Step 2: We compute the adjoints of $R_S$ and $U_S$. $R_S$ maps $C(F)$,
the continuous functions on $F$, to $C(S)$, the continuous functions
on $S$. So $R^*_S$ maps finite measures on $S$ to finite measures
on $F$. We have
$$\int f\, d(R_S^*\nu)=\int R_Sf\, d\nu=\int 1_S(x)f(x)\, \nu(dx),$$
and hence
\begin{equation} \label{Rstar}
R_S^*\nu(dx)=1_S(x)\, \nu(dx).
\end{equation}

$U_S$ maps $C(S)$ to $C(F)$, so $U_S^*$ maps finite measures on $F$
to finite measures on $S$. If $\nu$ is a finite measure on $F$,
then using \eqref{L213-E1}
\begin{align}
\int_S f\, d(U_S^*\nu)&=\int_F U_Sf\, d\nu=\int_F f\circ \vp_S(x)\, \nu(dx)
\label{L213-E2}\\
&=\int_F \Big( \sum_{T\in \sS_n} \frac{1_T(x)}{N_n(x)}\Big)f\circ \vp_S(x)\, 
\nu(dx)\nonumber\\
&=\sum_T \int_T \frac{f\circ \vp_S(x)}{N_n(x)}\, \nu(dx).\nonumber
\end{align}
Let $\vp_{T,S}:T\to S$ be defined to be 
the restriction of $\vp_S$ to $T$; this  is one-to-one and onto. If $\kappa$
is a measure on $T$, define its pull-back $\vp^*_{T,S}\kappa$ to be
the measure on $S$ given by
$$\int_S f\, d(\vp^*_{T,S}\kappa)=\int_T (f\circ \vp_{T,S})\, d\kappa.$$
Write
$$\nu_T(dx)=\frac{1_T(x)}{N_n(x)}\, \nu(dx).$$
Then \eqref{L213-E2} translates to
$$\int_S f \, d(U_S^*\nu)=\sum_T\int_T f\, \vp^*_{T,S}(\nu_T)(dx),$$
and thus
\begin{equation}\label{L213-E3}
U_S^*\nu=\sum_{T\in \sS_n} \vp^*_{T,S}(\nu_T).
\end{equation}

Step 3: We prove that if $\nu$ is a finite measure on $F$ such that
$\Theta^* \nu=\nu$ and $S\in \sS_n$, then
\begin{equation}\label{e-nu}
\nu(F)=\MF^n \int_S\frac{1}{N_n(x)}\, \nu(dx). 
\end{equation}
To see this, recall that $\vp^*_{T,V}(\nu_T)$ is a measure on $V$, and
then by \eqref{Rstar} and \eqref{L213-E3}
\begin{align*}
\Theta^* \nu&=\frac{1}{\MF^n} \sum_{V\in \sS_n}R_V^*U_V^*\nu\\
&=\frac{1}{\MF^n} \sum_{V\in \sS_n} \sum_{T\in \sS_n} \int 1_V(x)
\,\vp^*_{T,V}(\nu_T)(dx)\\
&=\frac{1}{\MF^n} \sum_V \sum_T \int \,\vp^*_{T,V}(\nu_T)(dx).
\end{align*}
On the other hand, using \eqref{L213-E1}
$$\nu(dx)=\sum_V \frac{1_V(x)}{N_n(x)}\, \nu(dx)=\sum_V \nu_V(dx).$$
Note that $\nu_V$ and 
$\MF^{-n} \sum_T \vp^*_{T,V}(\nu_T)$ are both supported on $V$, 
and the only way $\Theta^* \nu$ can equal $\nu$ is if
\begin{equation}\label{L213-E5}
\nu_V=\MF^{-n} \sum_{T\in \sS_n} \vp^*_{T,V}(\nu_T)
\end{equation} for each $V$.
Therefore   
\begin{align*}
\int_S \frac{1}{N_n(x)}\, \nu(dx)&=\nu_S(F)
= \MF^{-n}\sum_T \int 1_F(x)\, \vp^*_{T,S}(\nu_T)(dx)\\
&=\MF^{-n} \sum_T \int 1_F\circ \vp_{T,S}(x)\, \nu_T(dx)
=\MF^{-n}\sum_T \int \nu_T(dx)\\
&=\MF^{-n} \sum_T\int \frac{1_T(x)}{N_n(x)}\, \nu(dx)
=\MF^{-n} \int \nu(dx)=\MF^{-n} \nu(F).
\end{align*}
Multiplying both sides by $\MF^n$ gives \eqref{e-nu}.

Step 4: We show that if $\Theta f=f$, then 
\begin{equation}\label{e-GamQ}
\Theta^*(\Gamma(f,f))=\Gamma(f,f).
\end{equation}
Using Step 1, we have for $h\in C(F)\cap\sF$
\begin{align*}
\int_F h\, \Theta^*(\Gamma(f,f))(dx)&=\int_F \Theta h(x)\, \Gamma(f,f)(dx)
=2\sE(f,f\Theta h)-\sE(f^2, \Theta h)\\
&=2\sE(f, \Theta(fh))-\sE(\Theta f^2,h)=2\sE(\Theta f, fh)-\sE(f^2,h)\\
&=2\sE(f,fh)-\sE(f^2,h)=\int_F h\, \Gamma(f,f)(dx).
\end{align*}
This is the step where we used (b).

Step 5: We now prove (a). 
Note that if $g\in\sF\cap L^\infty (F)$ and 
$A=\{x\in F:g(x)=0\}$, then
$\Gam(g,g)(A)=0$ by Lemma~\ref{lem-zem}.  By applying this to the
function $g=f-U_SR_Sf$, which vanishes on $S$, 
and using the inequality 
\begin{align*}
\Big| \Gam(f,f)(B)^{1/2} -\Gam(U_SR_Sf,U_SR_Sf)(B)^{1/2}\Big|
 & \le\Gam(g,g)(B)^{1/2} \\
&\le\Gam(g,g)(S)^{1/2}=0,~~~\forall B\subset S,
\end{align*}
(see page 111 in \cite{FOT}), we see that 
\begin{equation}\label{e-Gam}
1_S(x)\Gam(f,f)(dx)=1_S(x)\Gam(U_SR_Sf,U_SR_Sf)(dx)
\end{equation}
for any $f\in\sF$ and $S\in\sS_n(F)$. 

Starting from  $U_TR_TU_SR_Sf=U_SR_Sf$, summing over $T\in \sS_n$ and dividing by
$\MF^n$ shows that $\Theta(U_SR_Sf)=U_SR_Sf$. Applying Step 4 with $f$ replaced
by $U_SR_Sf$,
$$\Theta^*(\Gamma(U_SR_Sf,U_SR_Sf))(dx)=\Gamma(U_SR_Sf, U_SR_Sf)(dx).$$
Applying Step 3 with $\nu=\Gamma(U_SR_Sf,U_SR_Sf)$, we see
\begin{align*}
 \sE(U_SR_Sf, U_SR_Sf) &=\Gamma(U_SR_Sf, U_SR_Sf)(F) \\
&=\MF^n \int_S \frac{1}{N_n(x)}\, \Gamma(U_SR_Sf, U_SR_Sf)(dx). 
\end{align*}
Dividing both sides by $\MF^n$, using the definition of $\sE^S$,
and \eqref{e-Gam},
\begin{equation}\label{ESint}
 \sE^S(R_Sf, R_Sf)=\int_S \frac{1}{N_n(x)}\, \Gamma(f, f)(dx).
\end{equation}
Summing over $S\in \sS_n$ and using \eqref{L213-E1} we obtain
$$ \sum_S \sE^S(R_Sf,R_Sf)=\int \Gamma(f,f)(dx)=\sE(f,f),$$
which is (a). 
\qed

\begin{corollary}\label{corGam}
If $\sE\in\fE$, $f\in\sF$, $S\in\sS_n(F)$, and $\Gam_S(R_Sf,R_Sf)$ is
the energy measure of $\sE^{S}$, then
$$\Gam_S(R_Sf,R_Sf)(dx)=\frac{1}{N_n(x)}\Gam(f,f)(dx), \q x \in S. $$ 
\end{corollary}

We finish this section with properties of sets of capacity zero for
$F$-invariant Dirichlet forms. 
Let  $A\subset F$ and $S\in \sS_n$. We define
\begin{equation}\label{defPOs} \Theta (A)=\vp_S^{-1}(\vp_S(A)). \end{equation}
Thus $\Theta (A)$ is the union of all the sets
that can be obtained from $A$ by local reflections. We can
check that $\Theta(A)$ does not depend on $S$, and that
$$ \Theta (A) = \{x:  \Theta (1_A)(x)>0 \} .$$

\begin{lemma}\label{lem-quasi}If  $\sE \in \fE$ then 
$$ \Cap(A) \le  \Cap(\Theta (A)) \le \MF^{2n} \Cap(A)$$ 
for all Borel sets $A\subset F$.  
\end{lemma} 

\proof The first inequality holds because we always have $A\subset
\Theta (A)$. To prove the second inequality it is enough to assume that $A$
is open since the definition of the capacity uses an infimum over open
covers of $A$, and $\Theta $ transforms an open cover of $A$ into an open
cover of $\Theta (A)$. If $u\in \sF$ and $u\geq 1$ on $A$, then
$\MF^n\Theta u\geq1$ on $\Theta (A)$. This implies the second inequality because
$\sE(\Theta u,\Theta u)\leq\sE(u,u)$, using that $\Theta $ 
is an orthogonal projection with
respect to $\sE$, that is, $\sE(\Theta f,g)=\sE(f,\Theta g)$.  \qed

\begin{corollary}\label{cor-quasi}
If $\sE \in \fE$, then $Cap(A)=0$ if and only if $Cap(\Theta (A))=0$. 
Moreover, if $f$ is quasi-continuous, then $\Theta f$ is quasi-continuous.
\end{corollary} 

\proof The first fact follows from Lemma~\ref{lem-quasi}. Then the
second fact holds because $\Theta $ preserves continuity of functions on
$\Theta $-invariant sets. \qed

\section{The Barlow-Bass and  Kusuoka-Zhou\\ Dirichlet forms}\label{BB-KZ}

In this section we prove that the Dirichlet forms associated with the
diffusions on $F$ constructed in \cite{BB1, BB4, KZ} are
$F$-invariant; in particular this shows that $\fE$ is non-empty and
proves Proposition \ref{ebbkz}.
A reader who is only interested in the uniqueness statement
in Theorem \ref{tmain} can skip this section.

\subsection{The Barlow-Bass processes}\label{SS:BB}

The constructions in \cite{BB1,BB4} were probabilistic and almost no
mention was made of Dirichlet forms.  Further,  
in \cite{BB4} the diffusion was
constructed on the unbounded fractal $\wt F$.  So before we can assert
that the Dirichlet forms are $F$-invariant, we need to discuss the
corresponding forms on $F$.  Recall the way the processes in
\cite{BB1,BB4} were constructed was to let $W^n_t$ be normally
reflecting Brownian motion on $F_n$, 
and to let
$X^n_t=W^n_{a_n t} $ for a suitable sequence $(a_n)$.
This sequence satisfied 
\begin{equation}\label{an-prop}
c_1 (\MF\rho_F/\LF^2)^n \le a_n \le c_2 (\MF\rho_F/\LF^2)^n,
\end{equation}
where $\rho_F$ is the resistance scale factor for $F$.
It was then
shown that the laws of the $X^n$ were tight and that resolvent
tightness held. Let $U^\lam_n$ be the $\lam$-resolvent operator for
$X^n$ on $F_n$. The two types of tightness were used to show there
exist subsequences $n_j$ such that $U^\lam_{n_j}f $ converges
uniformly on $F$ if $f$ is continuous on $F_0$ and that the $\bP^x$
law of $X^{n_j}$ converges weakly for each $x$. Any such a
subsequential limit point was then called a Brownian motion on the
GSC.  
The Dirichlet form for $W^n$ is $\int_{F_n} |\grad f|^2 \, d \mu_n$ 
and that for $X^n$ is 
$$\sE_n(f,f)= a_n \int_{F_n} |\grad f(x)|^2 \, \mu_n(dx),$$
both on $L^2(F, \mu_n)$.

Fix any subsequence $n_j$ such that the laws of the $X^{n_j}$'s converge, 
and the resolvents converge. If $X$ is the limit process and $T_t$
the semigroup for $X$, define
$$\sE_{BB}(f,f)=\sup_{t>0} \frac{1}{t} \langle f-T_tf,f\rangle$$
with the domain $\sF_{BB}$ 
being those $f\in L^2(F,\mu)$ for which the supremum is finite.

We will need the fact that if $U^\lam_n$ is the $\lam$-resolvent
operator for $X^n$ and $f$ is bounded on $F_0$, then $U^\lam_n f$
is equicontinuous on $F$. This is already known for the Brownian motion
constructed in \cite{BB4} on the unbounded fractal $\wt F$, but now 
we need it for
the process on $F$ with reflection on the boundaries of $F_0$. However the
proof is very similar to proofs in \cite{BB1,BB4}, so we will be brief. 
Fix  $x_0$ and suppose $x,y$ are in $ B(x_0,r)\cap F_n$. Then
\begin{align}\nn
U^\lam_n f(x)&=\bE^x \int_0^\infty e^{-\lam t} f(X^n_t)\, dt \\
&=\bE^x \int_0^{S^n_r} e^{-\lam t} f(X_t^n)\, dt +
\bE^x (e^{-\lam S^n_r}-1) U^\lam_n f(X^n_{S^n_r})
  +\bE^x U^\lam_n f(X^n_{S^n_r}), \label{Ucts}
\end{align} 
where $S^n_r$ is the time
of first exit from  $ B(x_0,r)\cap F_n$. The first term in 
\eqref{Ucts} is bounded by $\|f\|_\infty \bE^x S^n_r$.
The second term in \eqref{Ucts} is bounded by
$$\lam \| U^\lam_n f\|_\infty \bE^x S^n_r\leq \|f\|_\infty \bE^x S^n_r.$$ 
We have the same estimates in the case when $x$ is replaced by $y$, so
$$|U^\lam_n f(x)-U^\lam_n f(y)|\leq |\bE^x U^\lam_n f(X_{S^n_r}^n)
-\bE^y U^\lam_n f(X_{S^n_r}^n)|+\delta_n(r),$$
where $\delta_n(r)\to 0$ as $r\to 0$ uniformly in $n$
by \cite[Proposition 5.5]{BB4}.
 But
$z\to \bE^z U^\lam_n f(X_{S^n_r}^n)$ is harmonic in the ball of radius $r/2$
about $x_0$. Using the uniform elliptic Harnack inequality 
 for $X^n_t$ and the corresponding
uniform modulus of continuity for harmonic functions  (\cite[Section 4]{BB4}),
taking  $r=|x-y|^{1/2}$, and using the estimate for $\delta_n(r)$ gives
the equicontinuity. 

It is easy to derive from this that the limiting resolvent $U^\lam$
satisfies the property that $U^\lam f$ is continuous on $f$ whenever
$f$ is bounded.

\begin{theorem}\label{EBB} Each $\sE_{BB}$ is in $\fE$.
\end{theorem}

\proof
We suppose a suitable subsequence $n_j$ is fixed, and we write $ \sE$ for
the corresponding Dirichlet form $ \sE_{BB}$.
First of all, each $X^n$ is clearly conservative, so $T^n_t 1=1$. Since
we have $T_t^{n_j} f\to T_t f$ uniformly for each $f$ continuous, then
$T_t 1=1$. This shows $X$ is conservative, and $ \sE(1,1)=
\sup_t \langle 1-T_t1,1 \rangle =0$.

The regularity of 
$\sE$ follows from Lemma~\ref{lem-Feller} and the fact 
that the processes constructed in \cite{BB4} 
are $\mu$-symmetric Feller 
(see the above  discussion, \cite[Theorem~5.7]{BB4} and \cite[Section~6]{BB1}).
Since the process is a diffusion, the locality of $\sE$ follows from 
\cite[Theorem 4.5.1]{FOT}. 

The construction in \cite{BB1,BB4} gives a nondegenerate process, so
$ \sE$ is non-zero. Fix $\ell$ and let $S\in\sS_\ell(F)$. It is easy to
see from the above discussion that $U_SR_Sf\in\sF$ for any $f\in\sF$.
Before establishing the remaining properties of $F$-invariance, we
show that $\Theta _\ell$ and $T_t$ commute, where $\Theta _\ell$ is
defined in \eqref{defPO}, but with $\sS_n(F)$ replaced by
$\sS_\ell(F)$.  Let 
$\langle f,g\rangle_n$ denote $\int_{F_n} f(x)g(x)\, \mu_n(dx)$.  
The infinitesimal generator for $X^n$ is a constant times
the Laplacian, and it is clear that this commutes with $\Theta
_\ell$. Hence $U^\lam_n$ commutes with $\Theta _\ell$, or
\begin{equation}\label{Qcommute}
\langle \Theta _\ell U^\lam_n f,g\rangle_n
 =\langle U^\lam_n \Theta _\ell f,g \rangle_n.
\end{equation}
Suppose $f$ and $g$ are continuous and $f$ is nonnegative.  The left
hand side is $\langle U^\lam_n f, \Theta _\ell g \rangle_n$, and if
$n$ converges to infinity along the subsequence $n_j$, this converges
to 
$$\langle U^\lam f, \Theta _\ell g \rangle
 = \langle \Theta _\ell U^\lam f ,g\rangle.$$ 
The right hand side of \eqref{Qcommute}
converges to $\langle U^\lam \Theta _\ell f, g \rangle$ since $\Theta
_\ell f$ is continuous if $f$ is.  Since $X_t$ has continuous paths,
$t\to T_t f$ is continuous, and so by the uniqueness of the Laplace
transform, $\langle \Theta _\ell T_t f,g \rangle = \langle T_t\Theta
_\ell f, g\rangle$. Linearity and a limit argument allows us to extend
this equality to all $f\in L^2(F)$.  The implication
(c) $\Rightarrow$ (a) in Proposition~\ref{P2} implies that
$\sE\in \fE$.  \qed

\subsection{The Kusuoka-Zhou Dirichlet form}

Write $\sE_{KZ}$ for the Dirichlet form constructed in \cite{KZ}.
Note that this form is self-similar.

\begin{theorem}\label{EKZ} $\sE_{KZ} \in \fE$.
\end{theorem}

\proof One can see that $\sE_{KZ}$ satisfies Definition~\ref{deffE}
because of the self-similarity. The argument goes as follows.
Initially we consider $n=1$, and suppose $f\in\sF=\sD(\sEKZ)$. Then
\cite[Theorem 5.4]{KZ} implies $U_{S}R_S f \in \sF$ for any $S \in
\sS_1(F)$.  This gives us Definition~\ref{deffE}(1).

Let $S \in \sS_1(F)$ and $S=\Psi_i(F)$ where $\Psi_i$ is 
one of the contractions that define the self-similar structure on $F$, 
as in  \cite{KZ}. 
Then we have 
$$f\circ \Psi_i=(U_{S}R_S f)\circ \Psi_i=(U_{S}R_S f)\circ \Psi_j$$ 
for any $i,j$. Hence by \cite[Theorem 6.9]{KZ}, we have 
\begin{align*}
 \sEKZ(U_{S}R_Sf,U_{S}R_Sf) &= 
\rho_F \MF^{-1}\sum_{j}\sEKZ((U_{S}R_S f)\circ \Psi_j,(U_{S}R_S f)\circ \Psi_j)\\
 &= \rho_F \sEKZ(f\circ \Psi_i,f\circ \Psi_i).
\end{align*}
By \cite[Theorem 6.9]{KZ} this gives  Definition~\ref{deffE}(3), and 
moreover 
$$ \sE^{S}(f,f)=\rho_F \MF^{-1}\sEKZ(f\circ \Psi_i,f\circ \Psi_i).$$

Definition~\ref{deffE}(2) and the rest of  the conditions for $\sE_{KZ}$
to be in $\fE$  follow from (1), (3) and  the results of \cite{KZ}. 
The case  $n>1$ can be dealt with by using the self-similarity.    
\qed

\med {\bf Proof of Proposition \ref{ebbkz}}
This is immediate from Theorems  \ref{EBB} and \ref{EKZ}. 
\qed

\section { Diffusions associated with  $F$-invariant Dirichlet forms}
\label{Symm}

In this section we extensively use notation and definitions introduced in 
Section~\ref{SP}, especially Subsections~\ref{subGSC} and~\ref{subFinv}. 
We fix a Dirichlet form $\sE \in \fE$.  Let
$X=X^{(\sE)}$ be the associated diffusion, $T_t= T^{(\sE)}_t$ be the
semigroup of $X$ and $\bP^x=\bP^{x, (\sE)}$,  $x \in F -\sN_0$,
the associated probability laws. Here $\sN_0$ is a 
properly exceptional set for $X$. Ultimately (see Corollary
\ref{C1.3}) we will be able to define $\bP^x$ for all
$x \in F$, so that $\sN_0=\emptyset$.

\subsection{Reflected processes and the Markov property}

\begin{theorem}\label{T1}
Let $S \in \sS_n(F)$ and $Z=\vp_S(X)$.  Then $Z$
is a $\mu_S$-symmetric Markov process with Dirichlet form
$(\sE^{S}, \sF^{S})$, and semigroup $T^Z_t f = R_S T_t U_S f$.
Write $\wbP^y$ for the laws of $Z$; these are defined for
$y \in S-\sN^Z_2$, where $\sN^Z_2$ is a 
properly exceptional set for $Z$.
There exists a properly exceptional set $\sN_2$ for $X$ 
such that for any Borel set $A \subset F$, 
\begin{equation}\label{ZTrans}
 \wbP^{\vp_S(x)} (Z_t\in A) =\bP^x(X_t\in \vp_S^{-1}(A)), \q x \in 
F - \sN_2.
\end{equation}
\end{theorem}

\proof Denote $\vp=\vp_S$. 
We begin by proving that there exists a  properly exceptional
set $\sN_2$ for $X$ such that 
\begin{equation}\label{Lost1}
\bP^x(X_t\in
\vp^{-1}(A))=T_t1_{\vp^{-1}(A)}(x)=T_t1_{\vp^{-1}(A)}(y)=\bP^y(X_t\in
\vp^{-1}(A))
\end{equation}
whenever $A\subset S$ is Borel, $\vp(x)=\vp(y)$, and 
$x, y \in F -\sN_2$. 
It is sufficient to prove \eqref{Lost1} for a countable
base $(A_m)$ of the Borel $\sigma$-field on $F$.
Let $f_m = 1_{A_m}$. Since $T_t 1_{\vp^{-1}(A_m)}=T_t U_S f_m$, 
it is enough to prove that there exists a properly exceptional
set $\sN_2$ such that for $m \in \bN$,
\begin{equation}\label{sashasemigr}
 T_t U_S f_m(x)= T_t U_Sf_m(y), \q \hbox{ if $x,y \in F-\sN_2$ and
 $\vp(x)=\vp(y)$} . 
\end{equation}
By \eqref{URproj}, $\Theta (U_S f)=U_Sf$. Using Proposition~\ref{P2},
$$ \Theta T_t U_S f = T_t \Theta U_S f_m = T_t U_S f, $$
for $f \in L^2$, where the equality holds in the $L^2$ sense.

Recall that we always consider quasi-continuous modifications of
functions in $\sF$. By Corollary~\ref{cor-quasi}, 
$\Theta T_t U_S f_m$ is quasi-continuous. 
Since \cite[Lemma 2.1.4]{FOT} tells us that if two quasi-continuous functions 
coincide $\mu$-a.e., then they coincide q.e., we have that 
$\Theta (T_tU_Sf_m)=T_tU_Sf_m$ q.e. The definition of $\Theta$ implies
that $\Theta (T_t U_Sf_m)(x)=\Theta (T_tU_Sf_m)(y)$ whenever $\vp(x)=\vp(y)$,
so there exists a properly exceptional set $\sN_{2,m}$ such that
\eqref{sashasemigr} holds. Taking $\sN_2= \cup_m \sN_{2,m}$ gives
\eqref{Lost1}. Using Theorem 10.13 of \cite{Dyn}, 
$Z$ is Markov and has semigroup $T_t^Zf =R_ST_t(U_Sf)$. We take
$\sN^Z_2 =\vp(\sN_2)$.

Using \eqref{sashasemigr},
$U_S R_S T_t U_S f = T_t U_S f$,
and then
  $$\langle T^Z_tf,g \rangle_{S} = \langle R_S T_t U_S f, g \rangle_{S} 
= \MF^{-n} \langle U_S R_S T_t U_S f, U_S g \rangle 
= \MF^{-n} \langle T_t U_S f, U_S g \rangle.$$ 
This equals $\MF^{-n} \langle U_Sf, T_tU_S g\rangle$, and reversing
the above calculation, we deduce that  $\langle f, T_t^Z g\rangle
= \MF^{-n} \langle  U_S f, T_t U_S g \rangle,$
proving that $Z$ is $\mu_S$-symmetric.

To identify the Dirichlet form of $Z$ we note that 
$$ t^{-1} \langle T^Z_t  f -  f, f\rangle_{S}
 = \MF^{-n} t^{-1} \langle T_t U_S f -  U_S f, U_S f \rangle. $$
Taking the limit as 
$t\to 0$, and using \cite[Lemma 1.3.4]{FOT}, it follows that
$Z$ has Dirichlet form 
 $$\sE_Z(f,f)=\MF^{-n}\sE(U_S f,U_S f)=\sE^{S}(f,f).$$
\qed

\begin{lemma} \label{Zsymm}
Let $S, S' \in \sS_n$, $Z = \vp_S(X)$, and
$\Phi$ be an isometry of $S$ onto $S'$. Then if $x\in S-\sN$,
$$ \bP^x( \Phi(Z) \in \cdot) = \bP^{\Phi(x)}(Z \in \cdot). $$
\end{lemma}

\proof By Theorem~\ref{T1} and Definition~\ref{deffE}(2) $Z$ and
$\Phi(Z)$ have the same Dirichlet form. The result is then immediate
from \cite[Theorem 4.2.7]{FOT}, which states that two Hunt processes
are equivalent if they have the same Dirichlet forms, provided we
exclude an $F$-invariant set of capacity zero.
\qed

We say $S, S'\in \sS_n(F)$ are {\em adjacent} if
$S=Q \cap F$, $S'=Q'\cap F$ for $Q,Q'\in \sQ_n(F)$, and $Q\cap Q'$ is  
a $(d-1)$-dimensional set.
In this situation, let $H$ be the hyperplane separating $S, S'$.
For any hyperplane $H \subset \bR^d$, let $g_H : \bR^d \rightarrow \bR^d$ 
be reflection in $H$.  Recall the definition of $\pd_r D$, where $D$ 
is a finite union of elements of $\sS_n$.

\begin{lemma}\label{lem235}
  Let $S_1, S_2\in \sS_n(F)$ be adjacent,
let $D= S_1 \cup S_2$, let $B=\partial_r(S_1\cup S_2)$, and 
let $H$ be the hyperplane separating $S_1$ and $S_2$. 
Then there exists a properly exceptional set $\sN$ such that
if $x \in H \cap D -\sN$,
the processes $(X_t, 0\le t\le T_B)$ and 
$(g_H(X_t), 0\le t\le T_B)$ have the same law under $\bP^x$.
\end{lemma} 

\proof Let $f\in\sF$ with support in the interior of $D$. Then
Definition~\ref{deffE}(3) and Proposition~\ref{P4} imply that
$\sE(f,f)= \sE^{S_1}(R_{S_1} f, R_{S_1} f) +\sE^{S_2}(R_{S_2} f,
R_{S_2} f)$. Definition~\ref{deffE}(2) implies that $\sE(f,f)=
\sE(f\circ g_H, f\circ g_H)$.  Hence $(g_H(X_t), 0\le t\le T_B)$ has
the same Dirichlet form as $(X_t, 0\le t\le T_B)$, and so they have
the same law by \cite[Theorem 4.2.7]{FOT} if we exclude an
$F$-invariant set of capacity zero.  \qed

\subsection{Moves by $Z$ and $X$ }

At this point we have proved that the Markov process $X$ associated 
with the Dirichlet form $\sE \in \fE$ has strong symmetry 
properties. We now use these to obtain various global properties
of $X$. The key idea, as in \cite{BB4}, is to prove that certain
`moves' of the process in $F$ have probabilities which can be
bounded below by constants depending only on the dimension $d$.

We need a considerable amount of extra technical notation, based on
that in \cite{BB4}, which will only be used in this subsection.

We begin by looking 
at the process $Z=\vp_S(X)$ for some $S \in \sS_n$, where $n \ge 0$. 
Since our initial arguments are scale invariant, we can simplify 
our notation by taking $n=0$ and $S=F$ in the next definition. 

\begin{definition} \label{defLM} 
{\rm Let $1\le i,j \le d$, with $i \neq j$, and set
\begin{align*}
 H_i(t) &= \{x = (x_1, \dots, x_d): x_i =t \}, \, t \in \bR; \\ 
L_i     &= H_i (0) \cap [0, 1 /2]^d; \\
M_{ij}  &= \{ x \in [0,1]^d :\, x_i =0,\  \tfrac1 2 \leq x_j \leq 1, 
\mbox{ and } 0 \leq x_k \leq \tfrac12\mbox{ for } k \ne j
\}.
\end{align*}
}
\end{definition}

Let 
$$ \pd_e S = S \cap (\cup_{i=1}^d H_i(1)), \qq D= S - \pd_e S .  $$ 
We now define, for the process $Z$, the sets
$E_D$ and $Z_D$ as in \eqref{defED}. 
The next proposition says that the corners and slides of
\cite{BB4} hold for $Z$, provided that $Z_0 \in E_D$.

\begin{proposition}\label{Zmoves} There exists a constant
$q_0$, depending only on the dimension $d$, such that 
\begin{align} \label{m-corner}
 \wbP^x ( T^Z_{L_j} < \tau^Z_D ) &\ge q_0, \q x \in L_i \cap E_D, \\
\label{m-slide}
 \wbP^x ( T^Z_{M_{ij}} < \tau^Z_D ) &\ge q_0, \q x \in L_i \cap E_D. 
\end{align}These inequalities hold for any $n\ge0$ provided we modify 
Definition~\ref{defLM} appropriately. 
\end{proposition}

\proof Using Lemma \ref{Zsymm} this follows by the same reflection
arguments as those used in the proofs of Proposition 3.5 -- Lemma 3.10 
of \cite{BB4}. We remark that, inspecting these proofs, we can take  
$q_0 = 2^{-2d^2}$. \qed

\begin{figure}
\begin{center}
\begin{picture}(207,207)(0,0)\small\thinlines
\setlength{\unitlength}{69pt}
\put(0,0){\line(0,1){3}\line(1,0){3}}
\put(3,3){\line(0,-1){3}\line(-1,0){3}}
\put(1,1){\line(0,1){1}\line(1,0){1}}
\put(2,2){\line(0,-1){1}\line(-1,0){1}}
\linethickness{.25pt}
\multiput(0,2)(.025,0){80}{\line(0,1){1}}
\multiput(2,0)(.025,0){40}{\line(0,1){3}}
\linethickness{2pt}
\put(.52,1){\line(1,0){.43}}
\put(1,.53){\line(0,1){.43}}
\put(1,.47){\line(0,-1){.43}}
\put(1.04,.65){$A_0$}
\put(1.04,.2){$A_1$}
\put(.65,1.06){$A_1'$}
\put(1,1){\setlength{\unitlength}{1pt}\circle*{4}}
\put(1.03,.89){$v_*$}
\put(1,1){\rule{69pt}{69pt}}
\end{picture}
\end{center}
\caption{Illustration for Definition~\ref{move-n1} in the case of the
standard Sierpinski carpet and $n=1$. The complement of $D$ is shaded.
The half-face $A_1$ corresponds to a slide move, and the half-face $A_1'$
corresponds to a corner move. In this case $Q_*$ is the lower left
cube in $\sS_1$.}
\label{fig-move-n1}
\end{figure}
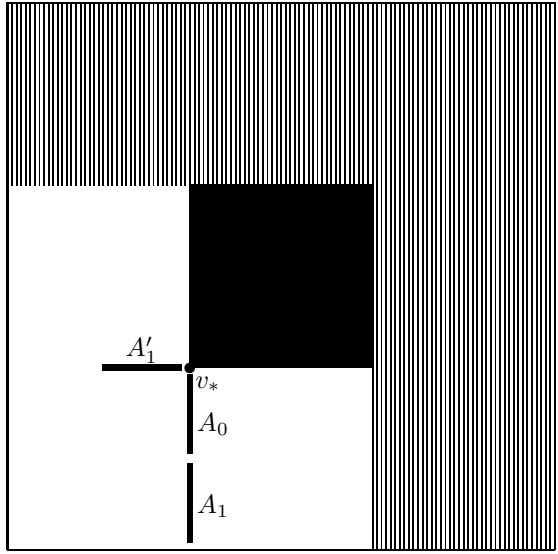

We now fix $n \ge 0$. 
We call a set $A \subset \bR^d$ a (level $n$) {\it half-face} if 
there exists $ i\in \{ 1,\ldots ,d\}$, $a=(a_1,\ldots , a_d)\in
\half  \bZ^d $ with $a_i\in \bZ$ such that
$$  A= \{ x: x_i=a_i \LF^{-n}, \quad 
a_j\LF^{-n} \le x_j \leq (a_j + {1/ 2})\LF^{-n} \quad\hbox{for }
j\ne i\} . $$(Note that a level $n$ half-face need not be a subset of $F$.)
For $A$ as above set $\iota (A)=i$. Let $\sA^{(n)} $ be the collection of
level $n$ half-faces, and 
$$ \sA^{(n)}_F =\{ A \in \sA^{(n)}: A \subset F_n \}. $$

We define a graph structure on $\sA^{(n)}_F$ by 
taking $\{A,B\}$ to be an edge if
$$ \dim (A\cap B) = d-2, \quad\hbox{and } A\cup B\subset Q 
\hbox{ for some } Q \in \sQ_n.  $$
Let $E(\sA^{(n)}_F)$ be the set of edges in $\sA^{(n)}_F$.
As in \cite[Lemma 3.12]{BB4} we have that the graph
$\sA^{(n)}_F$ is connected.
We call an edge $\{A,B\}$ an $i-j$ {\it corner }
if $\iota (A)=i$, $\iota (B)=j$, and $i\ne j$ and call $\{A,B\}$ an $i-j$
{\it slide } if $\iota (A) = \iota (B)=i$, and
the line joining the centers of $A$ and $B$ is parallel to the $x_j$ axis.
Any edge is either a corner or a slide; 
note that the move $(L_i,L_j)$ is 
an $i-j$ corner, while $(L_i, M_{ij})$ is an $i-j$ slide. 

For the next few results we need some further notation.

\sm \begin{definition}\label{move-n1}  {\rm 
Let $(A_0,A_1)$ be an edge in $E(\sA^{(n)}_F)$,
and $Q_*$ be a cube in $\sQ_n(F)$ such that $A_0 \cup A_1 \subset Q_*$. 
Let $v_*$ be the unique vertex of $Q_*$ such that $v_* \in A_0$,
and let $R$ be the union of the $2^d$
cubes in $\sQ_n$ containing $v_*$. 
Then there exist distinct  $S_i\in \sS_n$, $1\le i \le m$ such 
that $F \cap R = \cup_{i=1}^m S_i$.
Let $D = F \cap R^o$; thus 
$$ \ol D = F \cap R =  \cup_{i=1}^m S_i. $$ 
Let $S_*$ be any one of the $S_i$, and set $Z= \vp_{S_*}(X)$.
Write 
\begin{equation}
 \tau = \tau_D^X= \inf\{ t\ge 0: X_t \not\in D \} 
 = \inf\{t: Z_t \in \pd_r R \}. 
\end{equation}
Let
\begin{equation}
 E_D = \{ x \in D: \bP^x( \tau< \infty)=1 \}.
\end{equation}
} \end{definition}

We wish to obtain a lower bound for
\begin{equation}\label{XMlb}
  \inf_{x \in A_0 \cap E_D} \bP^x( T^X_{A_1} \le \tau). 
\end{equation}
By Proposition \ref{Zmoves} we have
\begin{equation}\label{Zmlb}
  \inf_{y \in A_0 \cap E_D} \wt \bP^y( T^Z_{A_1} \le \tau) \ge q_0. 
\end{equation}
$Z$ hits $A_1$ if and only if $X$ hits $\Theta(A_1)$, and one wishes 
to use symmetry to prove that, if $x \in A_0 \cap E_D$ then for some $q_1>0$
\begin{equation}\label{XZcomp}
  \bP^x( T^X_{A_1}  \le \tau) \ge q_1 
 \wt  \bP^x( T^Z_{A_1}  \le \tau) \ge q_1 q_0 .
\end{equation}
This was proved in \cite{BB4} in the context of reflecting Brownian
motion on $F_{n+k}$, but the proof used the fact that sets of 
dimension $d-2$ were polar for this process. 
Here we need to handle the possibility that there may be times $t$ such
that $X_t$ is in more than two of the $S_i$.
We therefore need to consider the
way that $X$ leaves points $y$ which are in several $S_i$. 

\begin{definition}\label{move-n2} {\rm
Let $y \in E_D$ be in exactly $k$ of the $S_i$,
where $1 \le k \le m$. 
Let $S_1', \dots, S'_k$ be the  elements of $\sS_n$ containing $y$.
(We do not necessarily have that 
$S_1$ is one of the $S'_j$.) 
Let $D(y) = \inter_r (\cup_{i=1}^k S'_i)$; so 
that $\ol{ D(y)} = \cup_{i=1}^k S'_i$.
Let $D_1$, $D_2$ be open sets in $F$ such that 
$y \in D_2 \subset \ol D_2 \subset D_1 \subset \ol D_1 \subset D(y)$.
Assume further that $\Theta(D_i)\cap D(y) = D_i$ for $i=1,2$, 
and note that we always have $\Theta(D_i) \supset D_i$. 
For $f \in \sF$ define
\begin{equation}
 \Theta^{D_1} f =  k^{-1} \MF^n 1_{D_1} \Theta f;
\end{equation}
the normalization factor is chosen so that
$\Theta^{D_1} 1_{D_1} = 1_{D_1}$.
} \end{definition}

\begin{figure}
\begin{center}
\begin{picture}(210,210)(0,0)\small\thinlines
\setlength{\unitlength}{70pt}
\put(0,0){\line(0,1){3}\line(1,0){3}}
\put(3,3){\line(0,-1){3}\line(-1,0){3}}
\put(1,1){\line(0,1){1}\line(1,0){1}}
\put(2,2){\line(0,-1){1}\line(-1,0){1}}
\linethickness{.25pt}
\multiput(0,2)(.025,0){80}{\line(0,1){1}}
\multiput(2,0)(.025,0){40}{\line(0,1){3}}
\linethickness{2pt}
\qbezier[40](0,0)(0,1)(0,2)
\qbezier[40](1,0)(1,1)(1,2)
\qbezier[20](0,2)(0.5,2)(1,2)
\qbezier[20](0,0)(0.5,0)(1,0)
\put(.5,1){
\qbezier[20](.4,0)(.4,0.7)(0,.7)
\qbezier[16](.3,0)(.3,0.4)(0,.4)
\qbezier[20](-.4,0)(-.4,0.7)(-0,.7)
\qbezier[16](-.3,0)(-.3,0.4)(-0,.4)
\qbezier[20](.4,0)(.4,-0.7)(0,-.7)
\qbezier[16](.3,0)(.3,-0.4)(0,-.4)
\qbezier[20](-.4,0)(-.4,-0.7)(0,-.7)
\qbezier[16](-.3,0)(-.3,-0.4)(0,-.4)
}
\put(.35,1){\setlength{\unitlength}{1pt}\circle*{4}}
\put(.39,1){$y$}
\put(1,1){\setlength{\unitlength}{1pt}\circle*{4}}
\put(1.03,.89){$v_*$}
\put(1,1){\rule{70pt}{70pt}}
{\put(.4,.1){\scriptsize$ D(y)$}}
{\put(.4,.42){\scriptsize$ D_1$}}
{\put(.4,.67){\scriptsize$ D_2$}}
\end{picture}
\end{center}
\caption{Illustration for  Definition~\ref{move-n2} in the case of
 the standard Sierpinski carpet and $n=1$. The complement of $D$
 is shaded, and the dotted lines outline $D(y)\supset D_1 \supset D_2$. 
}\label{fig-move-n2}
\end{figure}
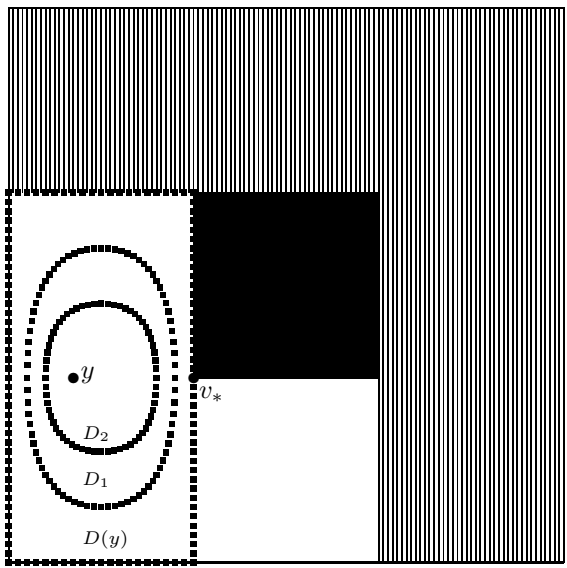

As before we define 
$ \sF_{D_1}\subset\sF$  as the closure of the set of functions $ \{ f \in \sF: \supp(f) \subset D_1 \} $. 
We denote by
$\sE_{D_1}$ the associated Dirichlet form and by $T_t^{D_1}$ the
associated semigroup, which are the Dirichlet form and the semigroup
of the process $ X$ killed on exiting $D_1$, by Theorems 4.4.3 and
A.2.10 in \cite{FOT}.  For convenience, we state the next lemma in the
situation of Definition \ref{move-n2}, although it holds under
somewhat more general conditions.

\begin{lemma}\label{lem-hitsym} 
Let $D_1$, $D_2$ be as in Definition \ref{move-n2}. \\
(a) Let $f \in \sF_{D_1}$. Then $\Theta^{D_1} f \in \sF_{D_1}$. 
Moreover, for all $f,g\in\sF_{D_1}$ we have 
$$\sE_{D_1}(\Theta^{D_1} f,g)=\sE_{D_1}(f,\Theta^{D_1} g)$$ 
and $T_t^{D_1}\Theta^{D_1} f=\Theta^{D_1} T_t^{D_1}f$. \\
(b) If $h \in \sF_{D_1}$ is 
harmonic (in the Dirichlet form sense) in $D_2$ then $\Theta^{D_1} h$ is
 harmonic (in the Dirichlet form sense) in $D_2$.  \\
(c) If $u$ is caloric in  $ D_2$, in the sense of Proposition~\ref{p-caloric}, 
then $\Th^{D_1} u$ is also caloric in  $D_2$.
\end{lemma}

\proof (a) By Definition \ref{deffE}, $\Theta f \in \sF$. Let $\psi$ be a
function in $\sF$ which has support in $D(y)$ and is 1 on $D_1$; such
a function exists because $\sE$ is regular and Markov.
Then $\psi \Theta f \in \sF$, and 
$\psi \Theta f = k\MF^{-n} \Theta^{D_1} f$. 
The rest of the proof follows from 
Proposition~\ref{P2}(b,c) because 
$\sE(\Theta^{D_1} f, g) = k^{-1} \MF^{n} \sE(\Theta f , g)$. \\
(b) Let $g \in \sF$ with $\supp(g) \subset D_2$. Then 
\begin{equation}\label{harm2}
 \sE(\Theta^{D_1} h, g) =  k^{-1} \MF^{n} \sE(\Theta h , g) 
=  k^{-1} \MF^{n} \sE(h, \Theta g)  = \sE(h, \Theta^{D_1} g)=0. 
\end{equation}
The final equality holds because $h$ is harmonic on $D_2$
and $\Theta^{D_1} g$ has support in $D_2$.
Relation \eqref{harm2} implies that $\Theta^{D_1} h$ is 
harmonic in $D_2$ by Proposition~\ref{s2p2}. 
\\ (c) We denote by $\ol T_t$ the semigroup of the process $\ol X_t$, 
which is $X_t$ killed at exiting $D_2$. 
The same reasoning as in (a) implies that 
$\ol T_t\Theta^{D_1} =\Theta^{D_1} \ol T_t$. Hence (c) follows 
 from  (a), (b) and Proposition~\ref{p-caloric}. \qed

Recall  from \eqref{Ndef} the definition of the ``cube counting'' 
function $N_n(z)$.
Define the related ``weight'' function  
$$ r_S(z) = 1_{S}(z) N_n(z)^{-1} $$ 
for each $S\in\sS_n(F)$. If no confusion can arise, we will denote
 $r_i(z)=r_{S'_i}(z)$. 

Let $(\sF^Z_t)$ be the filtration generated by $Z$. 
Since $\sF^Z_0$ contains all $\bP^x$ null sets, under the
law $\bP^x$ we have that $X_0=x$ is $\sF^Z_0$ measurable.

\begin{lemma}\label{equal-hit}
Let $y\in E_D$, $D_1$, $D_2$ be as in Definition \ref{move-n2}. 
Write $V = \tau^X_{D_2}$. \\
(a) If $U \subset \pd_F(D_2)$ satisfies $\Theta(U) \cap D(y) = U$, then
\begin{equation}\label{h-sym2}
\bE^y( r_i(X_{V}) 1_{(X_{V} \in U )} ) 
  = k^{-1} \wt \bP^{\vp_S(y)} ( Z_{V}  \in \vp_S(U)), \q 
    \hbox{ for $i=1,\ldots,k=N_n(y)$.}
\end{equation}
(b) For any bounded Borel function $f:D_1\to\bR$ and all $0\leq t\leq\infty$,
\begin{equation}\label{e-sym}
 \bE^y( f(X_{t\wedge V})| \sF^Z_{t \wedge V}) 
= \left(\Th^{D_1}f\right)(Z_{t \wedge V}).
\end{equation}  
In particular 
\begin{equation}\label{h-sym3}
 \bE^y (r_i(X_{t\wedge V})| \sF^Z_{t \wedge V}) = k^{-1}.
\end{equation}  
\end{lemma}

\proof 
Note that, by the symmetry of $D_2$, $V$ is a $(\sF^Z_t)$ stopping time.\\
(a) Let $f \in \sF_{D_1}$ be bounded,
and $h$ be the function with support in $D_1$ which equals $f$
in $D_1-D_2$, and is harmonic (in the Dirichlet form sense) inside $D_2$.
Then since $\vp_{S_i'}(y)=y$ for $1\le i \le k$, 
$$ \Theta^{D_1} h(y) = k^{-1} \sum_{i=1}^k h(\vp_{S_i'}(y)) = h(y). $$
Since $\Theta^{D_1} h$ is harmonic (in the Dirichlet form sense) 
in $D_2$ and since $y\in E_D$, 
we have, using Proposition~\ref{s2p2}, that 
$$  h(y) = \Theta^{D_1} h(y) = \bE^y ( \Theta^{D_1} h)(X_{V} )
= k^{-1} \bE^y \sum_{i=1}^k h(\vp_{S_i'}(X_{V} )). $$
Since $f=h$ on $\pd_F(D_2)$,
$$
 \bE^y(f(X_{V} ))= h(y)= k^{-1} \bE^y \sum_{i=1}^k f(\vp_{S_i'}(X_{V} )).
$$
Write $\delta_x$ for the unit measure at $x$, and define measures
$\nu_i(\om,dx)$ by
$$ \nu_1(dx) = \delta_{X_{V}}(dx), \q  \nu_2(dx) 
= k^{-1} \sum_{i=1}^k \delta_{\vp_{S_i'}(X_{V})}(dx)
= k^{-1} \sum_{i=1}^k \delta_{\vp_{S_i'}(Z_{V})}(dx) . $$
Then we have
$$ \bE^y \int f(x) \nu_1(dx) =  \bE^y \int f(x) \nu_2(dx) $$
for $f \in \sF_{D_1}$, and hence for all bounded Borel $f$ 
defined on $\pd_F(D_2)$. 
Taking $f= r_i(x) 1_U(x)$ 
then gives \eqref{h-sym2}. 

\noi (b) 
We can take the cube $S^*$ in Definition \ref{move-n1} to be $S'_1$.
If $g$ is defined on $S^*$ then $U_Sg$ is the unique extension of $g$
to $\ol{ D(y)}$ such
that $\Theta^{D_1} U_S g= U_S g$ on $\ol{ D(y)}$. Thus any function on
$S$ is the restriction of a function which is invariant with respect
to $\Theta^{D_1}$.  We will repeatedly use the fact that if 
$\Theta^{D_1} g = g $ then $g(X_t)=g(Z_t)$, and so also 
$g(X_{t\wedge V})=g(Z_{t\wedge V})$.

We break the proof into several steps. 

\noi Step 1.
Let $T_t^{D_2}$ denote the semigroup of $X$ stopped on exiting $D_2$, that is 
$$ T_t^{D_2}f(x)=\bE^x f(X_{t\wedge V}).$$ 
If  $f \in \sF_{D_1}$ is bounded, 
then Proposition~\ref{p-caloric} and Lemma~\ref{lem-hitsym} imply that q.e. in $D_2$ 
\begin{equation} \label{eTs}
T_t^{D_2}\Th^{D_1}f=\Th^{D_1} T_t^{D_2}f.
\end{equation}
Note that by Proposition~\ref{p-caloric} and \cite[Theorem 4.4.3(ii)]{FOT}, the 
notion  ``q.e.'' in $D_2$ coincides for the semigroups $T$, $T^{D_2}$ 
and $\ol  T$, where $\ol T$ is defined in Lemma~\ref{lem-hitsym}. 

\noi Step 2. 
If $f,g \in \sF_{D_1}$ are bounded and 
$ \Theta^{D_1} g = g $, then we have  $\Th^{D_1}(gf)
 =g\Th^{D_1}f$. Hence 
\begin{equation}\label{eTs-}
T_t^{D_2}(g\Th^{D_1}f)=T_t^{D_2}\Th^{D_1}(gf)=\Th^{D_1} T_t^{D_2}(gf).
\end{equation}

\noi Step 3. 
Let $\nu$ be a Borel probability measure on $D_2$.
Set $\nu^*= (\Th^{D_1})^*\nu$. Suppose that $\nu(\sN_2)=0$, where 
$\sN_2$ is defined in Theorem~\ref{T1}. 
If $f,g$ are as in the preceding paragraph, then we have  
\begin{align} \nn
 \bE^{\nu^*} g(Z_{t\wedge V})f(X_{t\wedge V}) 
&= \int_{D_2} T_t^{D_2} \big(gf\big)(x)\left(\Th^{D_1}\right)^*\nu(dx) \\
\nn
&= \int_{D_2}  \Th^{D_1}\big( T_t^{D_2}(gf)\big)(x) \nu(dx) \\
\nn
&= \int_{D_2}  T_t^{D_2}\big( g \Th^{D_1}f \big)(x) \nu(dx) \\
\nn 
&=  \bE^\nu g(Z_{t\wedge V})\Th^{D_1}f(X_{t\wedge V}) \\
\label{e:XZt}
&= \bE^\nu g(Z_{t\wedge V})\Th^{D_1}f(Z_{t\wedge V}),
\end{align}
where we use the definition of adjoint, \eqref{eTs-} to 
interchange $T^{D_2}$ and $\Th^{D_1}$, and that 
$g(X_{t\wedge V})=g(Z_{t\wedge V})$.

\noi Step 4. 
We prove by induction that if $\nu(\sN_2)=0$,
$m\ge 0$, $0<t_1<\cdots <t_m <t$, $g_1,\ldots,g_m$ 
are bounded Borel functions satisfying $ \Theta^{D_1} g_i = g_i $, 
and $f$ is bounded and Borel, then 
\begin{equation}\label{ind-m}
\bE^{\nu^*}(\Pi_{i=1}^{m}g_i(Z_{t_i\wedge V})) f(X_{t\wedge V}) = 
\bE^{\nu}(\Pi_{i=1}^{m}g_i(Z_{t_i\wedge V}))\Th^{D_1} f(Z_{t\wedge V}).
\end{equation}
The case $m=0$ is \eqref{e:XZt}. Suppose \eqref{ind-m} holds for
$m-1$. Then set
\begin{equation}\label{h-def}
 h(x) = \bE^{x}(\Pi_{i=2}^{m}g_i(Z_{(t_i-t_1)\wedge V})) f(X_{(t-t_1)\wedge V}).
\end{equation}
Write $\delta^*_x = (\delta_x)^*$. 
By \eqref{ind-m} for $m-1$, provided $x$ is such that $\delta_x^*(\sN_2)=0$,
\begin{align}\label{h-trans}
 \Theta^{D_1} h(x) 
 &=  \bE^{\delta_x^*} (\Pi_{i=2}^{m}g_i(Z_{(t_i-t_1)\wedge V})) f(X_{(t-t_1)\wedge V}) \\
 &= \bE^{x}(\Pi_{i=2}^{m}g_i(Z_{(t_i-t_1)\wedge V})) \Theta^{D_1} f(Z_{(t-t_1)\wedge V}).
\end{align}
So, using the Markov property, \eqref{e:XZt} and \eqref{h-trans}
\begin{align*}
\bE^{\nu^*}(\Pi_{i=1}^{m}g_i( & Z_{t_i\wedge V})) f(X_{t\wedge V})  
  = \bE^{\nu^*} g_1(Z_{t_1\wedge V}) h(X_{t_1\wedge V}) \\
&= \bE^{\nu} g_1(Z_{t_1\wedge V})  \Theta^{D_1} h(X_{t_1\wedge V}) \\
&= \bE^{\nu} g_1(Z_{t_1\wedge V}) \bE^{X_{t_1\wedge V}}\big(
   (\Pi_{i=2}^{m}g_i(Z_{(t_i-t_1)\wedge V})) \Theta^{D_1} f(Z_{(t-t_1)\wedge V})\big) \\
&= \bE^{\nu}(\Pi_{i=1}^{m}g_i(Z_{t_i\wedge V})) \Theta^{D_1} f(Z_{t\wedge V}),
\end{align*}
which proves \eqref{ind-m}. Therefore since $(\delta_x^*)^*=\delta_x^*$,
$$ \bE^{\delta_x^*} (\Pi_{i=1}^{m}g_i(Z_{t_i\wedge V})) f(X_{t\wedge V}) = 
\bE^{\delta_x^*} (\Pi_{i=1}^{m}g_i(Z_{t_i\wedge V}))\Th^{D_1}f(Z_{t\wedge V}), $$
and so 
$$ \bE^{\delta_x^*}( f(X_{t\wedge V})| \sF^Z_{t \wedge V}) 
 = \left(\Th^{D_1}f\right)(Z_{t \wedge V}). $$
To obtain \eqref{e-sym}, observe that $\delta_y^*=\delta_y$. 
Equation \eqref{h-sym3} follows 
since $\Th^{D_1}r_i(x)=k^{-1}$ for all $x\in D_1$. 
\qed

\def\twt{{t\wedge \tau}} 

\begin{corollary} Let $f: D(y) \to \bR$ be bounded Borel, and $t\ge 0$.
Then
\begin{equation}\label{e-sym2}
 \bE^y( f(X_{\twt})| \sF^Z_{\twt}) = 
 \left(\Th^{D(y)}f\right)(Z_{t \wedge \tau}).
\end{equation}  
\end{corollary}  

\proof This follows from Lemma \ref{equal-hit} by letting the regions 
$D_i$ in Definition \ref{move-n2} increase to $D(y)$. \qed

\sms 
Let $(A_0, A_1)$, $Z$  be as in Definition \ref{move-n1}. 
We now look at $X$ conditional on $\sF^Z$. 
Write $W_i(t) = \vp_{S_i}(Z_t)\in S_i$. For any $t$, we have 
that $X_\twt$ is at one of the points $W_i(\twt)$. Let
\begin{align*}
  J_i(t) &= \{ j: W_j(\twt) = W_i(\twt) \}, \\
   M_i(t) &= \sum_{j=1}^m 1_{( W_j(\twt)=W_i(\twt))}= \# J_i(t) , \\
  p_i(t) &= \bP^x( X_\twt = W_i(\twt) | \sF^Z_\twt) M_i(t)^{-1}=
 \bE^x( r_i(X_\twt) | \sF^Z_\twt) .
\end{align*}
Thus the conditional distribution of $X_t$ given
$\sF^Z_\twt$ is
\begin{equation}\label{eta-a}
\sum_{i=1}^k p_i(t)  \delta_{W_i(\twt)}.
\end{equation}
Note that by the definitions given above, we 
have $M_i(t)=N_n(W_i(t))$ for $0\le t< \tau$, 
which is the number of elements of $\sS_n$ that contain $W_i(t)$.

To describe the intuitive picture, we call the $W_i$ ``particles.'' 
Each $W_i(t)$ is a single point, and for
each $t$ we consider the collection of points 
$\{W_i(t), 1\le i \le m\}$. This is a finite set, but the number 
of distinct points depends on~$t$. In fact, we have 
$\{W_i(t), 1\le i \le m\}=\Theta\{X_t\}\cap D$. 
For each given~$t$, $X_t$ is equal to some of the $W_i(t)$.
If $X_t$ is in the $r$-interior of an element of $\sS_n$, 
then all the $W_i(t)$ are distinct, and so there are $m$
of them. In this case there is a single $i$ such that $X_t=W_i(t)$.
If $Z_t$ is in a lower dimensional face, then there can be fewer than~$m$
distinct points $W_i(t)$, because some of them coincide and we can have
$X_t=W_i(t)=W_j(t)$ for $i \neq j$. 
We call such a situation a ``collision.'' 
There may be many kinds of collisions 
because there may be many different lower dimensional faces that can be hit.

\begin{lemma} \label{pProp}
The processes $p_i(t)$ satisfy the following: \\
(a) If $T$ is any $(\sF^Z_t)$ stopping time satisfying $T \le \tau$ on
$\{ T < \infty\}$ then there exists $\delta(\omega)>0$ such that
\begin{equation*}
 p_i(T+h) =p_i(T) \q \hbox{ for } 0 \le h < \delta .
\end{equation*}
(b)  Let $T$ be any $(\sF^Z_t)$ stopping time satisfying $T \le \tau$ on
$\{ T < \infty\}$. Then for each $i=1, \dots k$,
\begin{equation*}
  p_i(T) = \lim_{s \to T-} M_i(T)^{-1} \sum_{j \in J_i(T)} p_j(s). 
\end{equation*}
\end{lemma}

\proof (a) 
Let $D(y)$ be as defined as in Definition \ref{move-n2},
and $D' = \vp_S( D(X_T))$.
Let 
$$ T_0= \inf\{ s \ge 0: Z_s \not\in D' \}, \q 
T_1= \inf\{ s \ge T: Z_s \not\in D' \}; $$
note that $T_1>T$ a.s.
Let $s>0$, $\xi_0$ be a bounded $\sF^Z_T$ measurable r.v., and
$ \xi_1 = \prod_{j=1}^m f_j( Z_{(T+t_j) \wedge T_1}), $
where $f_j$ are bounded and measurable, and 
$0\le t_1 < \dots < t_m \le s$.
Write $\xi'_1=  \prod_{j=1}^m f_j( Z_{(t_j) \wedge T_0})$.
To prove that $p_i( (T+s)\wedge T_1) = p_i(T)$ it is
enough to prove that
\begin{equation}\label{cexp-1}
 \bE^x \xi_0 \xi_1 r_i(X_{(T+s)\wedge T_1}) = \bE^x \xi_0 \xi_1 p_i(T). 
\end{equation}
However,
\begin{align}
\nonumber
 \bE^x \xi_0 \xi_1 r_i(X_{(T+s)\wedge T_1}) 
&=  \bE^x\Big( \xi_0 \bE( \xi_1  r_i(X_{(T+s)\wedge T_1}) | \sF^X_T)\Big) \\
\nonumber
&= \bE^x\Big( \xi_0 \bE^{X_T} ( \xi'_1  r_i(X_{s \wedge T_0}))\Big) \\
\label{cexp-3}
&= \bE^x\Big( \xi_0  \sum_{j}p_j(T) \bE^{W_j(T)}( \xi'_1  r_i(X_{s \wedge T_0}))\Big).
\end{align}
If $W_j(T) \not\in S_i$ then 
$$ \bE^{W_j(T)}( \xi'_1  r_i(X_{s \wedge T_0})) = 0. $$
Otherwise, by \eqref{h-sym3} we have
\begin{equation}
 \bE^{W_j(T)}( \xi'_1  r_i(X_{s \wedge T_0})) 
=  M_i(T)^{-1} \wt \bE^{Z_T} \xi_1'.  
\end{equation}
So,
\begin{align} \nn 
 \sum_{j}p_j(T)  \bE^{W_j(T)}( \xi'_1  r_i(X_{s \wedge T_0}))
&=  \sum_{j}p_j(T) 1_{(j \in J_i(T))} M_i(T)^{-1} \wt \bE^{Z_T} \xi_1' \\
\label{cexp-4}
&=  p_i(T) \wt \bE^{Z_T} \xi_1'.
\end{align}
Here we used the fact that $p_j(T)=p_i(T)$ if $j \in J_i(T)$.
Combining \eqref{cexp-3} and \eqref{cexp-4} we obtain 
\eqref{cexp-1}.  \\
(b) Note that 
$\sum_{ j \in J_i(T)} r_j(x)$ is constant   
in a neighborhood of $X_T$. Hence
$$ \lim_{s \to T-} \sum_{ j \in J_i(T)} r_j( X_{s})
 = \sum_{ j \in J_i(T)} r_j( X_T ), $$
and therefore
$$ \lim_{s \to T-} \sum_{ j \in J_i(T)} p_j(s)
 = \sum_{ j \in J_i(T)} p_j(T) = M_i(T) p_i(T), $$
where the final equality holds since $p_i(T)=p_j(T)$
if $W_i(T)=W_j(T)$. \qed

\begin{proposition} \label{Xmoves}
Let $(A_0, A_1)$, $Z$  be as 
in Definition \ref{move-n1}. 
There exists a constant $q_1>0$, depending only on $d$, such that
if $x \in A_0 \cap E_D$ and $T_0 \le \tau$ is a finite $(\sF^Z_t)$ 
stopping time, then
\begin{equation}\label{cond-move}
  \bP^x( X_{T_0} \in S | \sF^Z_{T_0}  ) \ge q_1.
\end{equation}
Hence 
\begin{equation}\label{eq-Xm}
  \bP^x( T^X_{A_1} \le \tau) \ge q_0 q_1. 
\end{equation}
\end{proposition}

\proof In this proof we restrict $t$ to $[0,\tau]$.  Lemma \ref{pProp}
implies that each process $p_i(\cdot)$ is a `pure jump' process, that
is it is constant except at the jump times.  (The lemma does not
exclude the possibility that these jump times might accumulate.)

Let
\begin{align*}
  K(t) &= \{ i: p_i(t) >0 \}, \\
  k(t) &= |K(t)|, \\ 
  p_{\min}(t) &= \min\{p_i(t): i \in K(t) \}= \min\{p_i(t):p_i(t)>0\}.
\end{align*}
Note that Lemma \ref{pProp} implies that if $p_i(t)>0$ then
we have $p_i(s)>0$ for all $s>t$. Thus $K$ and $k$ are non-decreasing
processes. Choose $I(t)$ to be the smallest $i$ such
that $p_{I(t)}(t)= p_{\min}(t)$.

To prove \eqref{cond-move} it is sufficient to prove that
\begin{equation}\label{eq-pk}
   p_{\min}(t) \ge 2^{-d k(t)} \ge 2^{-d 2^d}, \q  0 \le t \le \tau.
\end{equation}
This clearly holds for $t=0$, since $k(0)\ge 1$ and 
$p_i(0) = r_i(X_0)$, which is for each $i$ either zero or
at least $2^{-d}$.

Now let 
$$ T = \inf\{ t \le \tau:    p_{\min}(t) < 2^{-d k(t)} \}. $$
Since $p_i(T+h)=p_i(T)$ and $k(T+h)=k(T)$ for all sufficiently small 
$h>0$, we must have
\begin{equation}\label{eq-pmin2}
   p_{\min}(T) <  2^{-d k(T)}, \q \hbox { on $\{T < \infty$\}}. 
\end{equation}
Since $Z$ is a diffusion, $T$ is a predictable stopping time
so there exists an increasing sequence of stopping times
$T_n$ with $T_n<T$ for all $n$, and $T=\lim_n T_n$. 
By the definition of $T$, \eqref{eq-pk} holds for each $T_n$.
Let
$A =\{ \omega: k(T_n) < k(T)$ for all $n \}$.
On $A$ we have, writing $I=I(T)$, and using Lemma \ref{pProp}(b)
and the fact that $k(T_n) \le k(T)-1$ for all $n$,
\begin{align*}
    p_{\min}(T) &= p_I(T) = M_I(T)^{-1} \sum_{j \in J_I(T)} p_j(T) \\
 &= \lim_{ n \to \infty} M_I(T)^{-1} \sum_{j \in J_I(T)} p_j(T_n) 
\ge 2^{-d}  \lim_{ n \to \infty}  p_{\min}(T_n)  \\
&\ge  2^{-d}  \lim_{ n \to \infty} 2^{ - d k(T_n)} 
\ge  2^{-d} 2^{ - d (k(T)-1)} = 2^{- dk(T)}. 
\end{align*}
On $A^c$ we have
\begin{align*}
 p_{\min}(T) &= \lim_{ n \to \infty} M_I(T)^{-1} \sum_{j \in J_I(T)} p_j(T_n)  \\
&\ge  \lim_{ n \to \infty}  p_{\min}(T_n)  \\
&\ge  \lim_{ n \to \infty}  2^{ - d k(T_n)} =  2^{- dk(T)}. 
\end{align*}
So in both case we deduce that $p_{\min}(T)\ge  2^{- dk(T)}$,
contradicting \eqref{eq-pmin2}. It follows that $\bP(T < \infty)=0$,
and so \eqref{eq-pk} holds.

This gives \eqref{cond-move}, and using Proposition \ref{Zmoves} we then
obtain \eqref{eq-Xm}.  \qed

\subsection{Properties of $X$}\label{sec-moves} 

\begin{remark}\label{rem:doubling}
{\rm $\mu$ is a doubling measure, so for each Borel subset $H$ of $F$,
almost every point of $H$ is a point of density for $H$; see
\cite[Corollary IX.1.3]{torch}.  
}\end{remark}

\sms
Let $I$ be a face of $F_0$ and let $F'=F-I$.

\begin{proposition}\label{exitF}
There exists a set $\sN$ of capacity 0 such that if
$x\notin \sN$, then $\bP^x(\tau_{F'}<\infty)=1$.
\end{proposition}

\proof Let $A$ be the set of $x$ such that when the process starts at $x$,
it never leaves $x$.
Our first step is to show $F-A$ has positive measure. If not, for almost
every $x$, $T_tf(x)=f(x)$, so 
$$\frac{1}{t}\langle f-T_tf, f\rangle=0.$$
Taking the supremum over $t>0$, we have $\sE(f,f)=0$. This is true for
every $f\in L^2$, which contradicts $\sE$ being non-zero.

Recall the definition of $E_S$ in \eqref{defED}.
If $\mu(E_S\cap S)=0$  
for every $S \in \sS_n(F)$ and $n\ge 1$ then
$\mu(F-A)=0$. Therefore there must exist $n$ and $S\in \sS_n(F)$ such
that $\mu(E_S\cap S)>0$. 
Let $\eps>0$. By Remark \ref{rem:doubling}   
we can find $k\ge 1$ so that there exists $S'\in \sS_{n+k}(F)$ such that 
$$ \frac{\mu(E_S\cap S')}{\mu(S')}>1-\eps.$$
Let $S'' \in \sS_{n+k}$ be adjacent to $S'$ and contained in $S$, and
let $g$ be the  map that reflects $S'\cup S''$ across $S'\cap S''$.
Define
$$ J_i(S') = \cup\{ T: T \in \sS_{n+k+i}, \, T \subset \inter_r(S') \}, $$
and define $J_i(S'')$ analogously. We can choose $i$ large enough so that
\begin{equation}\label{muelb}
  {\mu(E_S\cap J_i(S'))} > (1- 2 \eps){\mu(S')}  .
\end{equation}

Let $x\in E_S\cap J_i(S')$. 
Since $x\in E_S$, the process
started from $x$ will leave $S'$ with probability one. 
We can find a finite sequence of moves (that is, corners or slides)
at level $n+k+i$ so that $X$ started at $x$ will exit $S'$
by hitting $S' \cap S''$. By Proposition \ref{Xmoves} 
the probability of $X$ following this sequence of moves is strictly 
positive, so we have
$$ \bP^x( X(\tau_{S'}) \in S' \cap S'' )>0. $$ 

Starting from $x\in E_S$, the process can never leave
$E_S$, so $X$ will leave $S'$ through $B=E_S\cap S'\cap S''$ with positive
probability.  
By symmetry, $X_t$ started from $g(x)$ will
leave $S''$ in $B$ with positive probability. So by  the strong Markov 
property, starting from $g(x)$, the process will leave $S$ with positive
probability. We conclude $g(x)\in E_S$ as well. Thus 
$g (E_S \cap J_i(S')) \subset E_S \cap J_i(S'')$, and so 
by \eqref{muelb} we have
\begin{equation*}
   {\mu(E_S\cap J_i(S''))} > (1- 2 \eps){\mu(S'')}  .
\end{equation*}
Iterating this argument, we have that for every $S_j \in \sS_{n+k}(F)$ 
with $S_j \subset S$,
\begin{equation*}
   \mu(E_S\cap S_j) \ge \mu(E_S \cap J_i(S_j)) \ge (1- 2 \eps)\mu(S_j). 
\end{equation*}
Summing over the $S_i$'s, we obtain
$$\mu(E_S\cap S)\geq (1-2\eps) \mu(S).$$
Since $\eps$ was arbitrary, then $\mu(E_S\cap S)=\mu(S)$. In other words,
starting from almost every point of $S$, the process will leave $S$.

By symmetry, this is also true for every element of $\sS_n(F)$
isomorphic to $S$.  Then using corners and slides (Proposition
\ref{Xmoves}), starting at almost any $x\in F$, there is positive
probability of exiting $F'$. We conclude that $E_{F'}$ has full
measure.

The function $1_{E_{F'}}$ is invariant so $T_t1_{E_{F'}}=1$, a.e. By
\cite[Lemma 2.1.4]{FOT}, $T_t(1-1_{E_{F'}})=0$, q.e. Let $\sN$ be the set
of $x$ where $T_t 1_{E_{F'}}(x)\ne 1$ for some rational $t$. If $x\notin \sN$, then
$\bP^x(X_t\in E_{F'})=1$ if $t$ is rational. By the Markov property, $x\in E_{F'}$.
\qed

\begin{lemma}\label{exits}
Let $U \subset F$ be open and non-empty.
Then $\bP^x(T_U < \infty)=1$, q.e.
\end{lemma} 

\proof This follows by Propositions \ref{Xmoves} and \ref{exitF}.
\qed

\subsection{Coupling}

\begin{lemma}\label{meastheo}
Let $(\Omega, \sF, \bP)$ be a probability space. Let $X$ and $Z$ be random variables 
taking values in  separable metric spaces $E_1$ and $E_2$, respectively,
each furnished with the Borel $\sigma$-field. 
Then there exists $F: E_2\times [0,1]\to E_1$
that is jointly measurable such that 
if $U$ is a random variable  whose distribution is uniform on $[0,1]$ 
which is independent of $Z$ 
and  $\wt X=F(Z,U)$, then $(X,Z)$
and $(\wt X,Z)$ have the same law.
\end{lemma}

\proof First let us suppose $E_1=E_2=[0,1]$. We will extend to the
general case later. Let $\bQ$ denote the rationals. For each $r\in [0,1]\cap
\bQ$, $\bP(X\leq r\mid Z)$ is a $\sigma(Z)$-measurable random variable,
hence there exists a Borel measurable function $h_r$ such that
$\bP(X\leq r\mid Z)=h_r(Z)$, a.s. For $r<s$ let
$A_{rs}=\{z: h_r(z)> h_s(z)\}$.
If $C=\cup_{r<s;\,r,s\in \bQ}
A_{rs}$, then $\bP(Z\in C)=0$. For $z\notin C$,
$h_r(z)$ is nondecreasing in $r$ for $r$ rational.
For $x\in [0,1]$,
define $g_x(z)$ to be equal to $x$ if $z\in C$ and equal to 
$\inf_{s>x, s\to x;\,s\in \bQ} h_s(z)$ otherwise.  For each $z$,
let $f_x(z)$ be
the right continuous  inverse to $g_x(z)$. Finally let
$F(z,x)=f_x(z)$.

We need to check that $(X,Z)$ and $(\wt X,Z)$ have the same distributions.
We have
\begin{align*}
\bP(X\leq x, Z\leq z)&=\bE[\bP(X\leq x\mid Z); Z\leq z]
=\lim_{s>x,s\in \bQ, s\to x} \bE[\bP(X\leq s\mid Z); Z\leq z]\\
&=\lim \bE[h_s(Z); Z\leq z]= \bE[g_x(Z); Z\leq z].
\end{align*}
On the other hand,
\begin{align*}
\bP(\wt X\leq x,Z\leq r)&=\bE[\bP(F(Z,U)\leq x\mid Z); Z\leq z]
=\bE[\bP(f_U(Z)\leq x\mid Z); Z\leq z]\\
&=\bE[\bP(U\leq g_x(Z)\mid Z); Z\leq z]=\bE[g_x(Z); Z\leq z].
\end{align*}

For general $E_1$, $E_2$, let $\psi_i$ be  bimeasurable
one-to-one maps from $E_i$ to $[0,1]$, $i=1,2$. Apply the above
to $\ol X=\psi_1(X)$ and $\ol Z=\psi_2(Z)$ to obtain a function
$\ol F$. Then $F(z,u)=\psi_1^{-1}\circ \ol F(\psi_2(z),u)$ will be  the
required function.
\qed

We say that $x,y\in F$ are {\sl $m$-associated}, and write
$x \tilm y$, if $\vp_S(x) =\vp_S(y)$ for some (and hence all) $S \in \sS_m$.
Note that by Lemma \ref{vp-prop} if $x \tilm y$ then also 
$x \sim_{m+1} y$. One can verify that this is the same as the definition 
of $x \tilm y$ given in \cite{BB4}.

\ms The coupling result we want is:

\begin{proposition}\label{basic-coup} 
(Cf.\ \cite[Theorem 3.14]{BB4}.)
Let $x_1, x_2 \in F$ with $x_1 \sim_n x_2$, where $x_1\in S_1\in \sS_n(F)$,
$x_2\in S_2\in \sS_n(F)$, 
and let $\Phi= \vp_{S_1}|_{S_2}$.
Then there exists a
probability space $(\Omega ,\sF ,  \bP)$ carrying 
processes $X_k$, $k=1,2$ and $Z$ with the following properties. \\
(a) Each $X_k$ is an $\sE$-diffusion started at $x_k$.\\
(b) $Z= \vp_{S_2}(X_2)=\Phi\circ \vp_{S_1}(X_1)$. \\
(c) $X_1$ and $X_2$ are conditionally independent given $Z$.
\end{proposition}

\proof 
Let $Y$ be the diffusion corresponding to
the Dirichlet form $\sE$ and let $Y_1, Y_2$ be processes 
such that $Y_i$ is equal in law to $Y$
started at $x_i$. 
Let $Z_1=\Phi\circ \vp_{S_1}(Y_1)$ and $Z_2=\vp_{S_2}(Y_2)$.
Since the Dirichlet
form for $\vp_{S_i}(Y)$ is $\sE^{S_i}$ and $Z_1, Z_2$ have the same starting point,
then $Z_1$ and $Z_2$ are equal in law.
Use Lemma \ref{meastheo} to find functions
$F_1$ and $F_2$ such that $(F_i(Z_i,U),Z_i)$ is equal in law to $(Y_i,Z_i)$, $i=1,2$,
if $U$ is an independent uniform random variable on $[0,1]$.

Now take a probability space supporting a process $Z$ with the same law as
$Z_i$ and two independent random variables $U_1, U_2$ independent of $Z$
which are uniform on $[0,1]$. Let $X_i=F_i(Z,U_i)$, $i=1,2$.
We proceed to show that the $X_i$ satisfy (a)-(c).

$X_i$ is equal in law to $F_i(Z_i,U_i)$, which is equal in law to $Y_i$, 
$i=1,2$, which establishes (a). Similarly $(X_i,Z)$ is equal  in
law to $(F(Z_i,U_i),Z_i)$, which is equal in law to $(Y_i,Z_i)$. Since
$Z_1=\Phi\circ \vp_{S_1}(Y_1)$ and $Z_2=\vp_{S_2}(Y_2)$,
it follows from the equality in law that $Z=\Phi\circ \vp_{S_1}(Y_1)$ and
$Z=\vp_{S_2}(Y_2)$.
This establishes (b). 

As $X_i=F_i(Z,U_i)$ for $i=1,2$, and $Z, U_1$, and $U_2$ are
independent, (c) is immediate.
\qed

Given a pair of $\sE$-diffusions $X_1(t)$ and $X_2(t)$
we define the coupling time 
\begin{equation}
  T_C(X_1, X_2) = \inf \{ t\ge 0: X_1(t)=X_2(t) \}.
\end{equation}

Given Propositions \ref{Xmoves} and \ref{basic-coup} we can now use
the same arguments as in \cite{BB4} to couple copies of $X$ started at
points $x, y \in F$, provided that $x \tilm y$ for some $m \ge 1 $.

\begin{theorem} \label{coupling}
Let $r>0$, $\eps>0$ and $r' = r/L_F^2$. 
There exist constants $q_3$ and $\delta$, depending only on the
GSC $F$, such that the following hold: \\
(a) Suppose $x_1$, $x_2 \in F$  with $||x_1-x_2||_\infty < r'$
and  $x_1 \tilm x_2$ for some $m \ge 1 $.
There exist $\sE$-diffusions $X_i(t)$, $i=1,2$, with
$X_i(0)=x_i$, such that, writing 
$$ \tau_i = \inf \{ t\ge 0: X_i(t) \not\in B(x_1, r) \}, $$
we have
\begin{equation}
  \bP\bigl( T_C(X_1, X_2) < \tau_1 \wedge \tau_2 \bigr) > q_3. 
\end{equation}
(b) If in addition $||x_1-x_2||_\infty < \delta r$ 
and  $x_1 \tilm x_2$ for some $m \ge 1 $ then 
\begin{equation}
   \bP\bigl( T_C(X_1, X_2) < \tau_1 \wedge \tau_2 \bigr) > 1 -\eps. 
\end{equation}
\end{theorem}

\proof Given  Propositions \ref{Xmoves} and \ref{basic-coup}, this
follows by the same arguments as in \cite{BB4}, p. 694--701. \qed

\subsection{Elliptic Harnack inequality}\label{SS-EHI}

As mentioned in Section $2.1$, there are two definitions of harmonic
that we can give.  We adopt the probabilistic one here.  Recall that a
function $h$ is harmonic in a relatively open subset $D$ of $ F$ if
$h(X_{t\land \tau_D'})$ is a martingale under $\bP^x$ for q.e.\ $x$
whenever $D'$ is a relatively open subset of $D$.

$X$ satisfies the {\sl elliptic Harnack inequality} if there exists a
constant $c_1$ such that the following holds: for any ball $B(x,R)$,
whenever $u$ is a non-negative harmonic function on $B(x,R)$ then
there is a quasi-continuous modification $\tilde u$ of $u$ that
satisfies
$$ \sup_{B(x,R/2)} {\tilde u} \le c_1 \inf_{B(x,R/2)} {\tilde u}.$$ 
We abbreviate ``elliptic Harnack inequality'' by ``EHI.'' 

\begin{lemma}\label{wk-holder}
Let $\sE$ be in $\fE$, $r \in (0,1)$, and $h$ 
be bounded and harmonic in $B=B(x_0,r)$.
Then there exists $\th>0$ such that
\begin{equation}\label{Hcont}
|h(x)-h(y)|\le  C \Big(\frac{|x-y|}{r}\Big)^\th (\sup_B |h|), \q x,y\in 
B(x_0, r/2), \q x \tilm y.
\end{equation}
\end{lemma}

\proof As in \cite[Proposition 4.1]{BB4}
this follows from the coupling in Theorem \ref{coupling} by
standard arguments. \qed

\begin{proposition}\label{HP1}
Let $\sE$ be in $\fE$ 
and $h$ be bounded and harmonic in $B(x_0,r)$. Then there exists
a set $\sN$ of $\sE$-capacity $0$ such that
\begin{equation}\label{stHcont}
|h(x)-h(y)|\le  C \Big(\frac{|x-y|}{r}\Big)^\th (\sup_B |h|), \q x,y\in 
B(x_0, r/2)-\sN. 
\end{equation}
\end{proposition}

\proof Write $B=B(x_0,r)$, $B'=B(x_0,r/2)$. By Lusin's theorem,
there exist open sets $G_n\downarrow$ such that $\mu(G_n)
\downarrow 0$, and $h$ restricted to $G_n^c\cap B'$ is continuous.
We will first show that $h$ restricted to any $G_n^c$ satisfies (\ref{Hcont})
except when one or both  of $x,y$ is in $\sN_n$, a  set of measure 0.
If $G=\cap_n G_n$,
then $h$ on $G^c$ is H\"older continuous
outside of $\cup \sN_n$, which is a  set of measure 0.
Thus $h$ is H\"older continuous on all of $B'$
outside of a set $E$ of measure 0.

So fix $n$ and let $H=G_n^c$.
Let $x,y$ be points of density for $H$; recall Remark \ref{rem:doubling}. 
Let $S_x$ and $S_y$ be appropriate isometries of an element of   $\sS_k$
such that $x\in S_x$, $y\in S_y$, and $\mu(S_x\cap H)/\mu(S_x)\geq \frac23$  and
the same for $S_y$. Let $\Phi$ be the isometry taking $S_x$ to $S_y$.
Then the measure of $\Phi(S_x\cap H)$ must be at least two thirds
the measure of $S_y$ and we already know the measure of $S_y\cap H$
is at least two thirds that of $S_y$.
Hence the measure of $(S_y\cap H)\cap (\Phi(S_x\cap H))$
is at least one third the measure of $S_y$.
So there must exist  points
$x_k\in S_x\cap H$ and $y_k=\Phi(x_k)\in S_y\cap H$
that are $m$-associated for some $m$. The inequality (\ref{Hcont})
holds for each pair
$x_k, y_k$. We do this for each $k$  sufficiently large and get 
sequences $x_k\in H$ tending
to $x$ and $y_k\in H$ tending to $y$. Since $h$ restricted to $H$ is
continuous, (\ref{Hcont}) holds for our given $x$ and $y$.

We therefore know that $h$ is continuous a.e. on $B'$.
We now need to show the continuity q.e., without modifying the function
$h$.
Let $x,y$ be two points in $B'$ for which $h(X_{t\land \tau_B})$ is
a martingale under $\bP^x$ and $\bP^y$. The set of points $\sN$ where
this fails has $\sE$-capacity zero. Let $R=|x-y|< r$ and let $\eps>0$.
Since $\mu(E)=0$, then by \cite[Lemma 4.1.1]{FOT}, for each $t$,
$T_t1_E(x)=T_t(x,E)=0$ for $m$-a.e.\ $x$. $T_t 1_E$ is in
the domain of $\sE$, so by \cite[Lemma 2.1.4]{FOT}, $T_t1_E=0$, q.e. 
Enlarge $\sN$ to include the null sets where $T_t1_E\neq0$ for some $t$
rational.
Hence if $x,y\notin \sN$, then with probability one with
 respect to both $\bP^x$ and $\bP^y$, we have
$X_t\notin E$ for $t$ rational.
Choose balls $B_x, B_y$ with radii in $[R/4,R/3]$  and centered at  $x$ and $y$,
resp.,  such
that $\bP^x(X_{\tau_{B_x}}\in \sN)=\bP^y(X_{\tau_{B_y}}\in \sN)=0$.
By the continuity of paths, we can choose $t$ rational and  small enough that
$\bP^x(\sup_{s\leq t}|X_s-X_0|>R/4)<\eps$ and the same with $x$
replaced by $y$. Then
\begin{align*}
|h(x)-h(y)|&=|\bE^x h(X_{t\land \tau_{B_x}})-\bE^y h(X_{t\land \tau_{B_y}})| \\
&\leq |\bE^x [h(X_{t\land \tau_{B_x}}); t<\tau_{B_x}]
-\bE^y [h(X_{t\land \tau_{B_y}}); t<\tau_{B_y}]|\, +2\eps\|{h}\|_\infty\\
&\leq C\Big(\frac Rr\Big)^\th\|{h}\|_\infty +4\eps\|{h}\|_\infty.
\end{align*}
The last inequality     above holds because  we
have $\bP^x(X_t\in \sN)=0$ and similarly  for $\bP^y$, points in $B_x$ are
 at most $2R$ from points in $B_y$, and
$X_{t\land \tau_{B_x}}$ and $X_{t\land \tau_{B_y}}$ are not in $E$ almost surely. 
Since $\eps$ is arbitrary, this
shows that except for $x,y$ in a set of capacity 0, we have (\ref{Hcont}).
\qed

\begin{lemma}\label{sm-ball}
Let $\sE \in \fE$. Then there exist constants $\kappa>0$, $C_i$, depending
only on $F$, such that if $0<r<1$, $x_0 \in F$, $y,z \in B(x_0,C_1 r)$ then
for all $0< \delta < C_1$,
\begin{equation}
 \bP^y( T_{B(z,\delta r)} < \tau_{B(x_0,r)}) > \delta^{\kappa}. 
\end{equation}
\end{lemma}

\proof This follows by using corner and slide moves, as in
\cite[Corollary 3.24]{BB4}. \qed

\begin{proposition}\label{HP2}
EHI holds for $\sE$, with constants depending only on $F$.
\end{proposition}

\proof Given Proposition \ref{HP1} and Lemma \ref{sm-ball}
this follows by the same argument as \cite[Theorem 4.3]{BB4}.
\qed

\begin{corollary} \label{irred}
(a) $\sE$ is irreducible.  \\
(b) If $\sE(f,f)=0$ then $f$ is a.e. constant.
\end{corollary}

\proof (a) If $A$ is an invariant set, then $T_t1_A=1_A$, or $1_A$ is
harmonic on $F$.  By EHI, either $1_A$ is never 0 except for a set of
capacity 0 or else it is 0, q.e. Hence $\mu(A)$ is either 0 or 1. So
$\sE$ is irreducible. \\
(b) The equivalence of (a) and (b) in this setting 
is well known to experts.
Suppose that $f$ is a function such that $\sE(f,f)=0$, and
that $f$ is not a.e. constant. Then using the contraction property
and scaling we can assume that $0\le f \le 1$ and there exist
$0<a<b<1$ such that the sets $A=\{x: f(x)<a \}$ and $B=\{x: f(x)>b
\}$ both have positive measure. Let $g= b \wedge (a \vee f)$; then
$\sE(g,g)=0$ also. By Lemma 1.3.4 of \cite{FOT}, for any $t>0$,
\begin{equation*}
 \sE^{(t)}(g,g) = t^{-1} \langle g - T_t g,g \rangle =0.
\end{equation*}
 So $ \langle g,T_t g \rangle = \langle g,g \rangle$.
By the semigroup property, $T_t^2=T_{2t}$, and hence $\langle T_t g,
T_t g \rangle =\langle g, T_{2t} g\rangle=\langle g,g \rangle$, from
which it follows that
$\langle g-T_t g ,g-T_t g \rangle =0$. This implies that $g(x) =
\bE^x g(X_t)$ a.e. Hence the sets $A$ and $B$ are invariant for
$(T_t)$, which contradicts the irreducibility of $\sE$. \qed

Given a Dirichlet form $(\sE, \sF)$ on $F$ we define
the {\em effective resistance} between subsets $A_1$ and $A_2$
of $F$ by:
\begin{equation}\label{def-reff}
 \Reff(A_1, A_2)^{-1} = \inf\{\sE(f,f): f\in \sF, f\mid_{A_1}=0, 
f\mid_{A_2}=1  \}.
\end{equation}
Let
\begin{equation}\label{defAt}
A(t)=\{ x \in F: x_1=t\}, \qq t\in [0,1].
\end{equation}

For $\sE \in \fE$ we set
\begin{equation}\label{def-enorm}
 ||\sE|| = \Reff(A(0), A(1))^{-1}. 
\end{equation}
Let $\fE_1 = \{ \sE \in \fE: ||\sE||=1 \}$.

\begin{lemma}\label{enorm}
If $\sE \in \fE$ then $||\sE||>0$. 
\end{lemma}

\proof 
Write $\sH$ for the set of functions $u$ on $F$ such that $u=i$ on 
$A(i)$, $i=0,1$. 
First, observe that $\sF\cap\sH$ is not empty. This is because,
by the regularity of $\sE$, there is a continuous function $u\in\sF$
such that $u\le0$ on the face $A(0)$ and $u\ge1$ on the opposite face
$A(1)$. Then the Markov property for Dirichlet forms says $0\vee
(u\wedge 1)\in\sF\cap\sH$.

Second, observe that by Proposition~\ref{exitF} and the symmetry,
$T_{A(0)}<\infty$ a.s., which implies that $(\sE,\sF_{A(0)})$ is a
transient Dirichlet form (see Lemma~1.6.5 and Theorem~1.6.2 in
~\cite{FOT}). Here as usual we denote 
$\sF_{A(0)}=\{f\in \sF: f|_{A(0)}=0\}$. 
Hence $\sF_{A(0)}$ is a Hilbert space with the norm
$\sE$. Let $u\in\sF\cap\sH$ and $h$ be its orthogonal projection onto
the orthogonal complement of $\sF_{A(0)\cup A(1)}$ in this Hilbert
space.  It is easy to see that $\sE(h,h)=||\sE||$.

If we suppose that $||\sE||=0$, then $h=0$ by
Corollary~\ref{irred}. By our definition, $h$ is harmonic in the
complement of $A(0)\cup A(1)$ in the Dirichlet form sense, and so by
Proposition~\ref{s2p2} $h$ is harmonic in the probabilistic sense and
$h(x) = \bP^x(X_{T_{A(0) \cup A(1)}} \in A(1))$. Thus, by the
symmetries of $F$, the fact that $h=0$ contradicts the fact that
$T_{A(1)}<\infty$ by Proposition~\ref{exitF}.

An alternative proof of this lemma starts with defining $h$
probabilistically and uses \cite[Corollary 1.7]{Ch} to show
$h\in\sF_{A(0)}$.  \qed

\subsection{Resistance estimates}\label{SS-Res}

Let now $\sE \in \fE_1$.
Let $S\in \sS_n$  and let $\gamma_n=\gamma_n(\sE)$ 
be the conductance across $S$.
That is, if $S=Q\cap F$ for $Q\in \sQ_n(F)$ and 
$Q=\{a_i\leq x_i\leq b_i, i=1, \ldots, d\}$, then
$$ \gamma_n=\inf\{\sE^S(u,u): u\in \sF^{S}, 
u\mid_{\{x_1=a_1\}}=0, u\mid_{\{x_1=b_1\}}=1\}.$$
Note that $\gam_n$ does not depend on $S$, and that
$\gamma_0=1$. Write $v_n=v_n^\sE$ for the minimizing
function. We remark that from the
results in \cite{BB3, McG} we have
$$ C_1 \rho_F^n \le \gam_n(\sEBB) \le C_2  \rho_F^n. $$

\begin{proposition} \label{P3.9} Let  $\sE \in \fE_1$.
Then for $n, m \ge 0$
\begin{equation}\label{gamineq}
 \gamma_{n+m}(\sE) \ge C_1 \gamma_m(\sE) \rho_F^{n}.
\end{equation}
\end{proposition}

\proof We begin with the case $m=0$. As in \cite{BB3} we compare the
energy of $v_0$ with that of a function constructed from $v_n$ and
the minimizing function on a network where each cube side $\LF^{-n}$
is replaced by a diagonal crosswire.

Write $D_n$ for the network of diagonal crosswires, as in \cite{BB3,
McG}, obtained by joining each vertex of a cube $Q \in \sQ_n$ to a
vertex at the center of the cube by a wire of unit resistance. Let
$R_n^D$ be the resistance across two opposite faces of $F$ in this
network, and let $f_n$ be the minimizing potential function.

Fix a cube $Q \in \sQ_n$ and let $S=Q\cap F$. 
Let $x_i$, $i=1, \dots 2^d$, be its vertices,
and for each $i$ let $A_{ij}$, $j=1, \dots d$, be the faces containing
$x_i$. Let  $A'_{ij}$ be the face opposite to  $A_{ij}$. Let
$w_{ij}$ be the function, congruent to $v_n$, which is 1 on $A_{ij}$
and zero on  $A'_{ij}$. Set
$$ u_i = \min\{ w_{i1}, \dots w_{id}\}. $$
Note that $u_i(x_i)=1$, and $u_i=0$ on $\cup_j A'_{ij}$.
Then
$$ \sE(u_i,u_i) \le \sum_j \sE(w_{ij},w_{ij}) = d \gam_n. $$

Write $a_i =f(x_i)$, and $\ol a= 2^{-d} \sum_i a_i$.
Then the energy of $f_n$ in $S$ is
$$ \sE^S_D(f_n,f_n) = \sum_i (a_i-\ol a)^2. $$
Now define a function $g_S: S \to \bR$ by
$$ g_S(y) = \ol a + \sum_i (a_i -\ol a) u_i(y). $$
Then
$$ \sE^S(g_S,g_S) \le C \sE(u_1,u_1) \sum_i  (a_i -\ol a)^2
 \le C \gamma_n  \sE^S_D(f_n,f_n). $$

We can check from the definition of $g_S$ that
if two cubes $Q_1$, $Q_2$ have a common face $A$ and $S_i=Q_i\cap F$,  then
$g_{S_1}=g_{S_2}$ on $A$.
Now define $g: F \to \bR$ by taking $g(x)=g_S(x)$ for
$x \in S$.
Summing over $Q \in  \sQ_n(F)$ we deduce that
$\sE(g,g) \le C \gamma_n (R^D_n)^{-1}$.
However, the function $g$ is zero on one face of $F$, and 1 on the
opposite face. Therefore
$$ 1=\gam_0 = \sE(v_0,v_0) \le \sE(g,g) \le  C \gamma_n (R^D_n)^{-1}
 \le C \gamma_n \rho_F^{-n}, $$
which gives \eqref{gamineq} in the case $m=0$.

The proof when $m\ge 1$ is the same, except we work in a cube $S \in \sS_m$
and use subcubes of side $L_F^{-n-m}$.  \qed

\begin{lemma} \label{gamma-sm}
We have
\begin{equation}\label{gam-comp}
 C_1 \gamma_n \le \gamma_{n+1} \le C_2 \gamma_n.
\end{equation}
\end{lemma}

\proof
The left-hand inequality is immediate from \eqref{gamineq}.
To prove the right-hand one, let first $n=0$.
By Propositions \ref{Xmoves} and \ref{exitF},   we deduce that
$v_0 \ge C_3>0$ on  $A(\LF^{-1})$; recall the definition in \eqref{defAt}.
Let $w= (v_0 \wedge C_3)/C_3$.
Choose a cube $Q  \in \sQ_1(F_1)$ between the hyperplanes
$A_1(0)$ and $A_1(\LF^{-1})$; $A_1(t)$ is defined in \eqref{defAt}. Then
\begin{align*}
 \gam_1 = \sE^{F_1}(v_1,v_1) &\le \sE^{F_1}(w,w) \le \sE(w,w) \\
 &= C_3^{-2} \sE( v_0 \wedge C_3, v_0 \wedge C_3)
 \le C_3^{-2} \sE(v_0,v_0) = C_4 \gam_0.
\end{align*}
Again the case $n\ge 0$ is similar, except we work in a
cube $S \in \sS_n$.
\qed

\sms Note that \eqref{gamineq} and \eqref{gam-comp} only give a one-sided
comparison between $\gam_n(\sE)$ and $\gam_n(\sEBB)$; however this 
will turn out to be sufficient. 

\ms Set
$$ \alpha=\log \MF/\log \LF,\qq
\beta_0 = \log (\MF\rho_F) /\log \LF. $$
By \cite[Corollary 5.3]{BB4} we have $\beta_0 \ge 2$,
and so $\rho_F \MF \ge \LF^2$. 
Let
$$ H_0(r) = r^{\beta_0}. $$

We now define a `time scale function' $H$ for $\sE$. First
note that by \eqref{gamineq} we have, for $n\ge 0$, $k \ge 0$.
\begin{equation}
\frac{\gam_n \MF^n}{\gam_{n+k} \MF^{n+k}}
 \le C \rho_F^{-k} \MF^{-k}.
\end{equation}
Since $\rho_F \MF \ge \LF^2 >1$ there exists $k\ge 1$ such that 
\begin{equation}
 \gam_n \MF^n  < \gam_{n+k} \MF^{n+k}, \q n \ge 0. 
 \end{equation}
Fix this $k$, let
\begin{equation}
   H(\LF^{-nk}) = \gam_{nk}^{-1} \MF^{-nk}, \q n \ge 0, 
 \end{equation}
and define $H$ by linear interpolation on each interval
$(\LF^{-(n+1)k}, \LF^{-nk})$. Set also $H(0)=0$. 
We now summarize some properties of $H$.

\begin{lemma}\label{hprop} There exist constants $C_i$
and $\beta'$, depending only on $F$ such that the 
following hold. \\
(a) $H$ is strictly increasing and continuous on $[0,1]$. \\
(b) For any $n, m \ge 0$
\begin{equation}\label{hreln}
  H(\LF^{-nk-mk}) \le C_1 H(\LF^{-nk}) H_0(\LF^{-mk}).
\end{equation}
(c) For $n \ge 0$
\begin{equation}\label{hreln2}
  H(\LF^{-(n+1)k}) \le H(\LF^{-nk}) \le C_2 H(\LF^{-(n+1)k}).
\end{equation}
(d) 
\begin{equation}\label{ftg}
   C_3 (t/s)^{\beta_0} \le \frac{H(t)}{H(s)} \le
  C_4 (t/s)^{\beta'}
\hbox { for } \,\, 0 < s \le t \le 1.
\end{equation}
In particular $H$ satisfies the `fast time growth' condition of
\cite {GT3} and \cite[Assumption 1.2]{BBKT}. \\
(e) $H$ satisfies `time doubling': 
\begin{equation}\label{timed}
   H(2r) \le C_5 H(r) \hbox { for } \,\, 0\le r \le 1/2.
\end{equation}
(f) For $r \in [0,1]$,
$$ H(r) \le C_6 H_0(r) . $$
\end{lemma}

\proof (a), (b) and (c) are immediate from the definitions of
$H$ and $H_0$, \eqref{gamineq} and \eqref{gam-comp}.
For (d), using \eqref{hreln} we have 
\begin{equation*}
 \frac{ H(\LF^{-kn})} {H( \LF^{-kn-km})}
 \ge C_7  \frac{ H(\LF^{-kn})}{  H(\LF^{-kn}) H_0(\LF^{-km})}
 = C_7 \LF^{km\beta_0}
= C_7  \Big(\frac{ \LF^{-kn}}{\LF^{-kn-km}} \Big)^{\beta_0},
\end{equation*}
and interpolating using (c) gives the lower bound in \eqref{ftg}.
For the upper bound,  using \eqref{gam-comp}, 
\begin{equation}
 \frac{ H(\LF^{-kn})} {H( \LF^{-kn-km})}
 \le C_8^{km} = \LF^{km \beta'} = 
 \Big(\frac{ \LF^{-kn}}{\LF^{-kn-km}} \Big)^{\beta'},
\end{equation}
where $\beta'= \log C_8/\log \LF$, and again using 
 (c) gives  \eqref{ftg}.
(e) is immediate from (d). 
Taking $n=0$ in (\ref{hreln}) and using (c) gives (f).
 \qed

We say $ \sE$ satisfies the condition RES$(H, c_1, c_2)$ if
for all $x_0\in F$, $r \in (0,\LF^{-1})$,
$$ c_1   \frac{H(r)}{r^\al} \le \Reff(B(x_0,r), B(x_0,2r)^c)
\le c_2 \frac{H(r)}{r^\al}. \eqno(RES(H,c_1,c_2))$$

\begin{proposition} \label{RESH}
There exist constants $C_1$, $C_2$, depending only on $F$,
such that $\sE$ satisfies  RES$(H,C_1,C_2)$.
\end{proposition} 

\proof Let $k$ be the smallest integer so that
$\LF^{-k} \le \half d^{-1/2} R$.
Note that if $Q \in \sQ_k$ and $x,y \in Q$, then
$d(x,y) \le d^{1/2} \LF^{-k} \le \half R$.
Write $B_0=B(x_0,R)$ and $B_1= B(x_0,2R)^c$.

We begin with the upper bound. Let $S_0$ be a cube
in $\sQ_k$ containing $x_0$: then $S_0 \cap F \subset B$.
We can find a chain of cubes $S_0, S_1, \dots S_n$
such that $S_n \subset B_1$ and $S_i$ is adjacent to $S_{i+1}$ for
$i=0, \ldots, n-1$. Let $f$ be the harmonic function
in $F-(S_0\cup B_1)$ 
which is 1 on $S_0$ and 0 on $B_1$. 
Let $A_0= S_0 \cap S_1$, and $A_1$ be the opposite face of
$S_1$ to $A_0$.
Then using the lower bounds for slides and corner moves, we have
that there exists $C_1 \in (0,1)$ such that
$f \ge C_1$ on $A_1$. So $g= (f-C_1)_+/(1-C_1)$
satisfies $\sE^{S_1}(g,g) \ge \gamma_k$.
Hence
$$ \Reff(S_0, B_1)^{-1} = \sE(f,f) \ge \sE^{S_1}(f,f)
 \ge (1-C_1)^{-2} \gamma_k, $$
and by the monotonicity of resistance
$$  \Reff(B_0, B_1) \le \Reff(S_0, B_1)
 \le C_2 \gamma_k^{-1}, $$ 
which gives the upper bound in (RES($H,c_1,c_2$)).

Now let $n=k+1$ and let $S \in \sQ_n$.
Recall from Proposition \ref{P3.9} the definition
of the functions $v_n$, $w_{ij}$ and $u_i$. By the symmetry
of $v_n$ we have that $w_{ij} \ge \half$ on the half of
$S$ which is closer to $A_{ij}$, and therefore
$u_i(x) \ge \half $ if  $||x-x_i||_\infty \le \half \LF^{-n}$.

Now let $y \in \LF^{-n} \bZ^d \cap F$, and let $V(y)$ be the union
of the $2^d$ cubes in $\sQ_n$ containing $y$. By looking at
functions congruent to $2 u_i  \wedge 1$ in each of the cubes in
$V(y)$, we can construct a function $g_i$ such that $g_i=0$ on $F -
V(y)$, $g_i(z)=1$ for $z \in F$ with $||z-y||_\infty \le \half
\LF^{-n}$, and $\sE(g_i,g_i)\le C \gamma_n$. We now choose $y_1, \dots
y_m$ so that $B_0 \subset \cup_i V(y_i)$: clearly we can take $m \le
C_5$. Then if $h= 1 \wedge (\sum_i g_i)$, we have $h=1$ on $B_0$ and
$h=0$ on $B_1$. Thus
$$  \Reff(B_0, B_1)^{-1} \le \sE(h,h)
\le \sE\Big(\sum g_i, \sum g_i\Big) \le C_6 \gamma_n, $$
proving the lower bound. \qed

\subsection{Heat kernel estimates}\label{SS-HK}

We write $h$ for the inverse of $H$, and $V(x,r)= \mu(B(x,r))$.
We say that $p_t(x,y)$  satisfies HK$(H;\eta_1, \eta_2,c_0)$
if for $x,y \in F$, $0<t\le 1$,
\begin{align*}
   p_t(x,y)
 &\ge c_0^{-1} V(x,h(t))^{-1} \exp(- c_0 ({H(d(x,y))/t})^{\eta_1}), \\
   p_t(x,y)
 &\le c_0  V(x,h(t))^{-1} \exp(- c_0^{-1} ({H(d(x,y))/t})^{\eta_2}).
\end{align*}

The following equivalence is proved in \cite{GT3}. 
(See also \cite[Theorem 1.3, $(a)\Rightarrow (c)$]{BBKT} 
 for a detailed proof of $(a)\Rightarrow (b)$, 
which is adjusted to our current setting.) 

\begin{theorem}\label{supplemain} 
Let $H: [0,1]\to[0,\infty)$ be a 
strictly increasing function with $H(1)\in (0,\infty)$ 
that satisfies (\ref{timed}) and (\ref{ftg}). Then 
the following are equivalent:
\nl (a) $(\sE,\sF)$ satisfies $(VD)$, $(EHI)$ and $(RES(H, c_1, c_2))$ 
for some $c_1,c_2>0$. 
\nl (b) $(\sE,\sF)$ satisfies $HK(H;\eta_1, \eta_2,c_0)$ 
for some $\al,\eta_1, \eta_2,c_0>0$.  \\
Further the constants in each implication are effective.
\end{theorem}

By saying that the constants are `effective' we mean that if,
for example (a) holds, then the constants $\eta_i$, $c_0$ in (b)
depend only on the constants $c_i$ in (a), and the constants
in (VD), (EHI) and  (\ref{timed}) and (\ref{ftg}).

\begin{theorem} \label{hkH}
$X$ has a transition density $p_t(x,y)$
which satisfies HK$(H;\eta_1,\eta_2, C)$, where $\eta_1=1/(\beta_0-1)$, 
$\eta_2=1/(\beta'-1)$, and the constant $C$ depends only on $F$.
\end{theorem}

\proof 
This is immediate from Theorem \ref{supplemain},
and Propositions \ref{HP2} and \ref{RESH}. \qed

Let
\begin{align}
 J_r(f) &=  r^{-\alpha} \int_F\int_{B(x,r)}
  |f(x)-f(y)|^2 d\mu(x) d\mu(y),\nonumber\\
 N_{H}^r(f)  &= H(r)^{-1} J_r(f), \nonumber\\
 N_H(f)  &= \sup_{0<r\le 1} N_H^r(f),\nonumber\\
 W_H
 &= \{ f\in L^2(F,\mu):  N_H(f)< \infty \}.\label{defWH}
\end{align}

We now use Theorem 4.1 of \cite{KS}, which we rewrite
slightly for our context. 
(See also Theorem 1.4 of \cite{BBKT}, 
which is adjusted to our current setting.) Let $r_j = L^{-kj}$,
where $k$ is as in the definition of $H$.

\begin{theorem} \label{twosum}
Suppose $p_t$ satisfies  HK$(H,\eta_1, \eta_2,C_0)$, and 
$H$ satisfies  (\ref{timed}) and (\ref{ftg}).
Then
\begin{equation}\label{e1.7}
C_1 \sE(f,f) \le \limsup_{j \to \infty}N_H^{r_j}(f) \le N_H(f) \le  C_2 \sE(f,f)
\q \hbox{ for all } f \in W_H,
\end{equation}
where the constants $C_i$ depend only on
the constants in  (\ref{timed}) and (\ref{ftg}),
and in  $HK(H;\eta_1,\eta_2, C_0)$. Further,
\begin{equation}\label{fwh}
 \sF = W_H.
\end{equation}
\end{theorem}

\begin{theorem}\label{ehequiv}
Let $(\sE,\sF) \in \fE_1$.
\newline (a) There exist constants
$C_1,C_2>0$ such that for all $r\in [0,1]$,
\begin{equation}\label{handho}
  C_1 H_0(r) \le H(r) \le C_2 H_0(r).
\end{equation}
(b) $W_H = W_{H_0}$, and there exist constants
$C_3,C_4$ such that
\begin{equation}\label{e-nh}
C_3 N_{H_0}(f) \le   \sE(f,f) \le C_4 N_{H_0}(f)
\q \hbox{ for all } f \in W_H.
\end{equation}
(c) $\sF = W_{H_0}$.  
\end{theorem}

\proof (a) We have $H(r) \le C_2 H_0(r)$
by Lemma \ref{hprop}, and so
\begin{equation}\label{nheq1}
  N_H(f) \ge C_2^{-1} N_{H_0}(f).
\end{equation}
Recall that $(\sE_{BB},\sF_{BB})$ is (one of) the Dirichlet forms
constructed in \cite{BB4}.
By \eqref{nheq1} and \eqref{fwh} we have
$\sF \subset \sF_{BB}$.
In particular, the function $v_0^\sE \in \sF_{BB}$ \  (see Subsection~\ref{SS-Res}).

Now let
\begin{equation*}
  A = \limsup_{k \to \infty} \frac{H(r_k)}{ H_0(r_k)};
\end{equation*}
we have $A \le C_2$.

Let $f \in \sF$. Then by Theorem \ref{twosum}
\begin{align*}
 \sE_{BB}(f,f)
&\le C_3 \limsup_{j \to \infty} H_0(r_j)^{-1} J_{r_j} (f) \\
&= C_3 \limsup_{j \to \infty}
 \frac{H(r_j)}{ H_0(r_j)}  H(r_j)^{-1} J_{r_j} (f) \\
&\le C_3 \limsup_{j \to \infty} A  N^{r_j}_H(f) \le C_4 A \sE(f,f).
\end{align*}

Taking $f =v_0^\sE$,
\begin{equation}\label{}
 1 \le  \sE_{BB}( v_0^\sE,v_0^\sE ) \le  C_4 A
  \sE( v_0^\sE,v_0^\sE ) = C_4 A.
\end{equation}
Thus $A \ge C_5 = C_4^{-1}$.
By Lemma \ref{hprop}(c) we have, for $n, m \ge 0$,
$$ \frac{H(r_{n+m})}{H_0(r_{n+m})} \le C_6
 \frac{H(r_{n})}{H_0(r_{n})}. $$
So, for any $n$
$$  \frac{H(r_{n})}{H_0(r_{n})} \ge C_6^{-1} A \ge
  C_5/C_6, $$
and (a) follows.

\noindent (b) and (c) are then immediate by Theorem \ref{twosum}. \qed

\sm 
\begin{remark} {\rm
\eqref{handho} now implies that $p_t(x,y)$ satisfies HK$(H_0, \eta_1, \eta_1, C)$
with $\eta_1= 1/(\beta_0-1)$.
}
\end{remark}

\section {Uniqueness}\label{Uni}

\begin{definition}\label{defHilb}{\rm  Let $W =W_{H_0}$ be as
defined in \eqref{defWH}.
Let $\sA, \sB \in \fE$. We say
$\sA \le \sB$ if
$$ \sB(u,u)-\sA(u,u) \ge 0 \hbox { for all }  u \in W. $$
For  $\sA, \sB \in \fE$ define
$$ \sup(\sB| \sA) = \sup \left\{ \frac{\sB(f,f)}{\sA(f,f)}: f \in W\right\}, $$
$$ \inf(\sB| \sA) = \inf \left\{ \frac{\sB(f,f)}{\sA(f,f)}: f \in W\right\}, $$
$$ \h(\sA,\sB) = \log \left( \frac{\sup(\sB| \sA)}{ \inf(\sB| \sA)}\right); $$
$\h$ is Hilbert's projective metric and we have $\h(\theta\sA, \sB)=
\h(\sA,\sB)$ for any $\theta\in (0,\infty)$.
Note that $h(\sA,\sB)=0$ if and only if $\sA$ is a nonzero
constant multiple of $\sB$.
}
\end{definition}

\begin{theorem}\label{hmbound}
There exists a constant $C_F$, depending only on the GSC $F$,
such that if $\sA, \sB \in \fE$ then
$$  \h(\sA,\sB) \le C_F. $$
\end{theorem}

\proof Let $\sA'=\sA/||\sA||$,  $\sB'=\sB/||\sB||$. Then
$\h(\sA,\sB)=  \h(\sA',\sB')$. By Theorem \ref{ehequiv} there exist
$C_i$ depending only on $F$ such that \eqref{e-nh} holds for both
$\sA'$ and $\sB'$. Therefore
$$ \frac{\sB'(f,f)}{\sA'(f,f)} \le \frac{C_2}{C_1}, \qq
\hbox { for } f \in W, $$
and so $\sup(\sB'|\sA') \le C_2/C_1$. Similarly, $\inf(\sB'|\sA')
\ge C_1/C_2$, so $\h(\sA',\sB') \le 2 \log(C_2/C_1)$. \qed

\med {\bf Proof of Theorem \ref{tmain}}\label{proof-tmain}
By Proposition \ref{ebbkz} we have that $\fE$ is non-empty.
 
Let $\sA, \sB \in \fE$, and  $\lam= \inf(\sB|\sA)$. 
Let $\delta>0$ and $\sC=(1+\delta)\sB - \lam \sA$.
By Theorem \ref{dfgen},
 $\sC$ is a local regular Dirichlet form on $L^2(F,\mu)$ and $\sC\in \fE$.
Since 
$$\frac{\sC(f,f)}{\sA(f,f)}=(1+\delta)\frac{\sB(f,f)}{\sA(f,f)}-\lam, \qquad
f\in W,$$  we obtain
$$ \sup(\sC| \sA) = (1+\delta) \sup(\sB| \sA)-\lam, $$
and
$$ \inf(\sC| \sA) =  (1+\delta) \inf(\sB| \sA)-\lam = 
\delta  \lam. $$
Hence for any $\delta>0$,
$$ e^{\h(\sA, \sC)}
= \frac{ (1+\delta) \sup(\sB| \sA)-\lam}{\delta \lam}
 \ge \frac1\delta\left ( e^{\h(\sA,\sB)}-1\right). $$
If $h(\sA, \sB)>0$, this is not bounded as $\delta \to 0$, contradicting
Theorem \ref{hmbound}.  We must therefore have $h(\sA,\sB)=0$, which
proves our theorem.
\qed

\med {\bf Proof of Corollary \ref{C1.3}}
Note that Theorem \ref{tmain} implies
that the $\bP^x$ law of $X$ is uniquely defined, up to scalar
multiples of the time parameter, for all $x\notin \sN$, where $\sN$ is
a set of capacity 0.  If $f$ is continuous and $X$ is a Feller
process, the map $x\to \bE^x f(X_t)$ is uniquely defined for all $x$
by the continuity of $T_tf$. By a limit argument it is uniquely
defined if $f$ is bounded and measurable, and then by the Markov
property, we see that the finite dimensional distributions of $X$
under $\bP^x$ are uniquely determined. Since $X$ has continuous paths,
the law of $X$ under $\bP^x$ is determined.  (Recall that the the
processes constructed in \cite{BB4} are Feller processes.)
\qed

\begin{remark}\label{rem:chargeboundary} 
{\rm In addition to (H1)-(H4), assume that the 
$(d-1)$-dimensional fractal 
$F\cap\{x_1=0\}$ also satisfies the conditions corresponding
to (H1)-(H4).  (This assumption is used in \cite[Section 5.3]{HK}.). 
Then one can show $\Gamma (f,f)(F\cap \partial F_0)=0$
for all $f\in \sF$ where $\Gamma(f,f)$ is the energy measure for
$\sE\in \fE$ and $f\in \sF$.  Indeed, by the uniqueness we know
that $\sE$ is self-similar, so the results in \cite{HK} can be
applied.  For $h$ given in \cite[Proposition 3.8]{HK}, we
have $\Gamma (h,h)(F\cap\partial [0,1]^d)=0$ by taking $i\to\infty$
in the last inequality of \cite[Proposition 3.8]{HK}. For general
$f\in \sF$, take an approximating sequence $\{g_m\}\subset \sF$ as
in the proof of Theorem 2.5 of \cite{HK}. Using the inequality
\begin{align*}
|\Gamma(g_m,g_m)(A)^{1/2} - \Gamma(f,f)(A)^{1/2}| 
 &\le \Gamma(g_m-f,g_m-f)(A)^{1/2} \\
 &\le 2 \sE(g_m-f,g_m-f)^{1/2},
\end{align*}
(see page 111 in \cite{FOT}), we conclude that 
$\Gamma (f,f)(F\cap\partial [0,1]^d)=0$.
Using the self-similarity, we can also prove that the energy measure 
does not charge 
the image of $F\cap\partial [0,1]^d$ by any of the contraction maps.
}
\end{remark}

\begin{remark}\label{rem:convergence}
{\rm One question left over from \cite{BB1,BB4} is whether
the sequence of approximating reflecting Brownian motions
used to construct the Barlow-Bass processes converges. Let
$\wt X^n_t=X^n_{c_nt}$, where $X^n$ is defined in Subsection
\ref{SS:BB} and $c_n$ is a normalizing constant. We choose
$c_n$ so that the expected time for $\wt X^n$ started at 0
to reach one of the faces not containing 0 is one. There
will exist subsequences $\{n_j\}$ such that there is resolvent
convergence for $\{\wt X^{n_j}\}$  and also weak convergence,
starting at every point in $F$. Any of the subsequential limit points
will have a Dirichlet form that is a constant multiple of one of
the $\sE_{BB}$. By virtue of the normalization and our uniqueness
result, all the limit points are the same, and therefore
the whole sequence $\{\wt X^n\}$ converges, both in the sense of
resolvent convergence and in the sense of weak convergence for
each starting point. }
\end{remark}

\bigskip
\noindent {\bf Acknowledgment. } The authors thank 
Z.-Q.~Chen, 
M.~Fukushima, 
M.~Hino, 
V.~Metz, and 
M.~Takeda 
for valuable discussions, and 
D.~Croydon for correcting some typos.


\begin{thebibliography}{}

\bibitem{AO} Alexander, S. and Orbach, R.: 
Density of states on fractals: \lq\lq fractons\rq\rq. 
{\sl J. Physique (Paris) Lett.} {\bf 43} (1982), 625-631. 

\bibitem{bar} Barlow, M.T.: 
{\sl Diffusions on fractals}.  
Lectures on probability theory and statistics (Saint-Flour, 1995), 1--121, 
Lecture Notes in Math. {\bf 1690}, Springer, Berlin, 1998.

\bibitem {BB1} Barlow, M.T., Bass, R.F.:  The construction of Brownian
motion on the Sierpinski carpet.  {\sl Ann. Inst.  H. Poincar\'e \bf
25} (1989) 225--257.

\bibitem {BB3} Barlow, M.T., Bass, R.F.: On the resistance of the Sierpinski
carpet. {\sl Proc. R. Soc. London A.} {\bf 431} (1990) 345--360.

\bibitem {BB4} Barlow, M.T., Bass, R.F.: Brownian motion and harmonic
analysis on Sierpinski carpets.  {\sl Canad. J. Math.}
{\bf 54} (1999), 673--744.

\bibitem {BBK1}
Barlow, M.T., Bass, R.F., Kumagai, T.: Stability of parabolic
Harnack inequalities on metric measure spaces.
{\sl J. Math. Soc. Japan} (2) {\bf 58} (2006), 485--519.

\bibitem{BK} Barlow, M.T., Kumagai, T.:  
Random walk on the incipient infinite cluster on trees.  
{\sl Illinois J. Math.}, {\bf 50} (2006), 33--65. 

\bibitem{BJKS} 
Barlow, M.T., J\'arai, A.A., Kumagai, T., Slade, G.:  
Random walk on the incipient infinite cluster for oriented percolation in high 
dimensions. {\sl Commun. Math. Phys.}, {\bf 278} (2008), 
385--431. 

\bibitem {BBK2} Barlow, M.T., Bass, R.F., Kumagai, T.: 
Note on the equivalence of parabolic Harnack inequalities and
heat kernel estimates. \hfill\break
{\tt http://www.math.uconn.edu/$\sim$bass/papers/phidfapp.pdf}

\bibitem {BBKT} Barlow, M.T., Bass, R.F., Kumagai, T., Teplyaev, A.: 
Supplementary notes for \lq\lq Uniqueness of Brownian motion on
Sierpinski carpets".\hfill\break
 {\tt http://www.math.uconn.edu/$\sim$bass/papers/scuapp.pdf}

\bibitem{BP} 
Barlow, M.T., Perkins, E.A.: 
Brownian Motion on the Sierpinski Gasket. 
{\sl Probab. Theory Rel. Fields} {\bf 79} (1988), 543--623.

\bibitem{HBA} Ben-Avraham, D., Havlin, S.: 
{\sl Diffusion and reactions in fractals and disordered systems.}   
Cambridge University Press, Cambridge, 2000. 

\bibitem {BH} Bouleau, N., Hirsch, F.: 
{\sl Dirichlet forms and analysis on Wiener space.}
de Gruyter Studies in Mathematics, 14. Walter de Gruyter and Co.,
Berlin, 1991.

\bibitem {Ch} Chen, Z.-Q.: 
On reflected Dirichlet spaces. {\sl Probab. Theory Rel. Fields} {\bf
94} (1992), 135--162.

\bibitem {Chnew} Chen, Z.-Q.: 
On notions of harmonicity.
{\sl Proc. Amer. Math. Soc.}, to appear. 

\bibitem{Dyn} Dynkin, E.B.: {\sl Markov Processes - I.}
Springer, Berlin,  1965.

\bibitem {FOT} Fukushima, M., Oshima, Y., Takeda, M.: 
{\sl Dirichlet Forms and Symmetric Markov Processes.}
de Gruyter, Berlin, 1994.

\bibitem{FS} Fukushima, M., Shima, T.: 
On a spectral analysis for the  Sierpi\'nski gasket. 
{\sl Potential Anal.} {\bf 1} (1992), 1--35.

\bibitem{G} Goldstein, S.:  
Random walks and diffusions on fractals. 
Percolation theory and ergodic theory of infinite particle systems 
(Minneapolis, Minn., 1984--1985), 121--129, 
IMA Vol. Math. Appl. {\bf 8}, Springer, New York, 1987. 

\bibitem {GT3} Grigor'yan, A., Telcs, A.: 
Two-sided estimates of heat kernels in metric measure spaces.
In preparation. 

\bibitem {HMT} Hambly, B.M., Metz, V., Teplyaev, A.: 
 Admissible refinements of energy on finitely ramified fractals.
{\sl J. London Math. Soc.} {\bf 74} (2006), 93--112.

\bibitem {HK} Hino, M., Kumagai, T.:  
A trace theorem for Dirichlet forms on fractals. 
{\sl J. Func. Anal.} {\bf 238} (2006), 578--611.

\bibitem{Kes1}
Kesten, H.: Subdiffusive behavior of random walk on a random cluster. 
{\sl Ann. Inst. H. Poincare Probab. Statist.}, 
{\bf 22} (1986), 425Ð487. 

\bibitem{Kigsg}  Kigami, J.: 
A harmonic calculus on the Sierpinski space. 
{\sl Japan J.\ Appl.\ Math.} {\bf 6} (1989), 259--290.

\bibitem{Kigpcf} Kigami, J.: 
A harmonic calculus for p.c.f.~self-similar sets.
{\sl Trans. Amer. Math. Soc. } {\bf 335} (1993), 721--755.

\bibitem{kig} Kigami, J.: 
{\sl Analysis on fractals}. Cambridge Univ. Press, Cambridge, 2001.

\bibitem{KN}
Kozma, G., Nachmias, A.: 
The Alexander-Orbach conjecture holds in high dimensions. 
Preprint 2008.

\bibitem {KS}
Kumagai, T., Sturm, K.-T.: Construction of diffusion processes on
fractals, $d$-sets, and general metric measure spaces.
{\sl J. Math. Kyoto Univ.} {\bf  45}  (2005),  no. 2, 307--327.

\bibitem{Kus0} Kusuoka, S.: 
A diffusion process on a fractal. 
Probabilistic methods in mathematical physics (Katata/Kyoto, 1985), 
251--274, Academic Press, Boston, MA, 1987. 

\bibitem{KZ} Kusuoka, S., Zhou, X.Y.: 
Dirichlet forms on fractals: Poincar\'e constant and resistance.
{\sl Probab. Theory Related Fields \bf 93} (1992), no. 2, 169--196.

\bibitem{L} Lindstr\o m, T.:  
Brownian motion on nested fractals. 
{\sl Mem.\ Amer.\ Math.\ Soc.} {\bf 420} (1990), 1--128.

\bibitem{Man} Mandelbrot, B.: {\it The Fractal Geometry of Nature}. W.H. Freeman,
San Francisco, 1982.

\bibitem{Me1} Metz, V.:  
{Renormalization contracts on nested fractals.}
{\sl J. Reine Angew. Math.} \textbf{480} (1996), 161--175.

\bibitem{McG} McGillivray, I.:
Resistance in higher-dimensional Sierpi\'nski carpets.
{\sl Potential Anal.} {\bf 16} (2002),  no. 3, 289--303.

\bibitem{Osada2006} Osada, H.:  
{Singular time changes of diffusions on Sierpinski carpets}.
{\sl Stochastic Process. Appl.}   {\bf116} (2006),   675--689.

\bibitem{Peirone2000} Peirone, R.:  
Convergence and uniqueness problems for Dirichlet forms on 
fractals. 
{\sl Boll. Unione Mat. Ital. Sez. B Artic. Ric. Mat. (8)} 
{\bf 3} (2000),   431--460.

\bibitem{RT} Rammal, R., Toulouse G.: 
Random walks on fractal structures and percolation clusters.  
{\sl J. Physique Lettres} {\bf 44} (1983) L13--L22.

\bibitem{RS} Reed, M., Simon, B.: {\sl Methods of modern mathematical
physics. I. Functional analysis}.  Academic Press, 1980.

\bibitem{RW} Rogers, L.C.G., Williams, D.: 
{\sl Diffusions, Markov Processes, and Martingales. Volume one: Foundations},
2nd ed. Wiley, 1994.

\bibitem{Ru} Rudin, W.: {\sl Functional analysis}.  McGraw-Hill,  1991.

\bibitem{Sab1} Sabot, C.:  
Existence and uniqueness of diffusions on finitely ramified self-similar
fractals.
{\sl Ann. Sci. \'Ecole Norm. Sup. (4)}  {\bf30} (1997), 605--673.

\bibitem{Schum} Schmuland, B.:   
On the local property for positivity preserving coercive forms. 
Dirichlet forms and stochastic processes (Beijing, 1993), 345--354, 
de Gruyter, Berlin, 1995.

\bibitem{Stri} Strichartz, R.S.: 
{\sl Differential Equations on Fractals: a Tutorial.} 
Princeton University Press, Princeton, NJ, 2006. 

\bibitem{torch} Torchinsky, A.: 
{\sl Real-Variable Methods in Harmonic Analysis.}
Academic Press, Orlando FL, 1986.

\end{thebibliography}
\end{document}